\documentclass[12pt]{amsart}

\usepackage{amssymb}
\usepackage{a4wide}
\usepackage{graphicx}
\usepackage{tikz-cd}
\usepackage{caption}
\usepackage{subcaption}

\newcommand{\br}{{\rm br}}
\newcommand{\tr}{{\rm tr}}
\newcommand{\C}{\mathbb{C}}
\newcommand{\R}{\mathbb{R}}
\newcommand{\Q}{\mathbb{Q}}
\newcommand{\Z}{\mathbb{Z}}
\newcommand{\N}{\mathbb{N}}
\newcommand{\B}{\mathbb{B}}
\newcommand{\Br}{\mathcal{B}}

\newcommand{\CHL}{\mathcal{C}}

\renewcommand{\P}{\mathbb{P}}

\newtheorem{prop}{Proposition}
\newtheorem{thm}{Theorem}
\theoremstyle{remark}
\newtheorem{rk}{Remark}

\def\cqfd{\mbox{}\nolinebreak\hfill$\Box$\medbreak\par}

\usepackage{color}

\title{Volumes of 3-ball quotients as intersection numbers}

\author{Martin Deraux}

\date{May 25, 2019}

\begin{document}

\begin{abstract}
  We give an explicit description of the 3-ball quotients constructed
  by Couwenberg-Heckman-Looijenga, and deduce the value of their
  orbifold Euler characteristics. For each lattice, we also give a
  presentation in terms of generators and relations.
\end{abstract}

\maketitle

\section{Introduction}
Let $X=G/K$ be an irreducible symmetric space of non-compact type. It
is a well known fact originally due to Borel~\cite{borel} that $G$
contains lattices, i.e. discrete subgroups such that $\Gamma\backslash
X$ has finite volume. In fact, $\Gamma$ can be chosen so that
$\Gamma\backslash X$ is compact or non-compact.

The standard construction of lattices comes from arithmetic, as we
briefly recall. Take a linear algebraic group $H$ defined over $\Q$,
and denote by $H^0_\R$ the connected component of the identity in the
group of real points $H_\R$. Assume that there is a surjective
homomorphism $\varphi:H^0_\R\rightarrow G$ with compact kernel, and
consider the group $\Gamma=\varphi(H^0_\R\cap H_\Z)$. It is a standard
fact that $\Gamma$ is then a lattice in $G$ (this is essentially a
result by Borel and Harish-Chandra~\cite{bhc}).

By definition, a lattice $\Gamma'$ in $G$ is called arithmetic if
there exists $H,\varphi,\Gamma$ as above such that $\Gamma$ and
$\Gamma'$ are commensurable in the wide sense, i.e. possibly after
replacing $\Gamma$ by $g\Gamma g^{-1}$ for some $g\in G$, the
intersection $\Gamma\cap\Gamma'$ has finite index in both $\Gamma$ and
$\Gamma'$. It follows from important work of Margulis~\cite{margulis},
Corlette~\cite{corlette}, Gromov-Schoen~\cite{gromovschoen} that if
$X$ is not a real or complex hyperbolic space, then every lattice in
$G$ is actually arithmetic.

In the case $X=H^n_\R$, $G=PO(n,1)$, several constructions of
non-arithmetic lattices are known, but no general structure theory for
lattices has been worked out. It follows from a construction of Gromov
and Piatetski-Shapiro~\cite{gromovpiat} that there exist
non-arithmetic lattices in $G$ for arbitrary $n\geqslant 2$, and that
there are infinitely many commensurability classes in each dimension.

The case $X=H^n_\C$, $G=PU(n,1)$ is even further from being
understood. There is currently no generalization of the
Gromov-Piatetski-Shapiro construction to the complex hyperbolic case,
and in fact (for $n>2$) only finitely many commensurability classes of
non-arithmetic lattices are known, only in very low dimension; there
are currently 22 known classes in $PU(2,1)$, see~\cite{thealgo}, and 2
known classes in $PU(3,1)$ see~\cite{chl3d_1}.

The first examples were constructed by Mostow~\cite{mostowpacific},
generalized by Deligne-Mostow~\cite{delignemostow}, then the list was
expanded~\cite{dpp2},~\cite{thealgo},~\cite{chl3d_1}. Some recent
constructions rely on the use of fundamental domains (and heavy
computational machinery), but most examples have been given
alternative constructions using orbifold uniformization
(see~\cite{derauxklein},~\cite{derauxabelian}).

It turns out most known examples are in fact in a list of lattices
that was produced by Couwenberg, Heckman and Looijenga~\cite{chl}
(their list contains representatives of 17 out of the 22 classes in
$PU(2,1)$, and both classes in $PU(3,1)$). For the sake of brevity,
we refer to their lattices as CHL lattices, and to the corresponding
quotients as CHL ball quotients.

An explicit description of the quotient of all the 2-dimensional CHL
lattices can be obtained by combining the results
in~\cite{delignemostow} and~\cite{derauxklein}. The goal of this paper
is to give an explicit description of the quotient for all
3-dimensional CHL lattices. In principle a similar description can of
course be worked out for higher-dimensional examples (recall that CHL
lattices only exist in dimension at most 7).

Using this description, we compute orbifold Euler characteristics of
the 3-dimensional CHL ball quotients. Recall that the orbifold Euler
characteristic is a universal multiple of the volume, namely
$$
  Vol(\Gamma\backslash\B^n)=\frac{(-4\pi)^n}{(n+1)!}\chi^{orb}(\Gamma\backslash\B^n),
$$ 
if the metric is normalized to have holomorphic sectional curvature
$-1$ (this is an orbifold version of the Chern-Gauss-Bonnet formula,
see~\cite{satake}). 

Since most of these lattices are arithmetic, one could in principle
compute their covolumes by using the Prasad
formula~\cite{prasadvolumes} (for all but one lattice, namely the
non-arithmetic one). Note however that Prasad's formula gives the
covolume of a specific lattice in each commensurability class (the
so-called principal arithmetic lattices); unfortunately the relation
of a given lattice to the principal arithmetic lattice in its
commensurability class can be difficult to make explicit.
In fact, our volume computations should make it possible to relate
arithmetic CHL lattices to the corresponding principal arithmetic
groups in their commensurability class. It may also be useful in order
to distinguish commensurability classes of non-arithmetic lattices,
using the Margulis commensurator theorem and volume estimates, in the
spirit of the arguments in~\cite{thealgo}.

Note that volumes of Deligne-Mostow ball quotients (which are special
cases of CHL lattices) were already known. They were computed by
McMullen~\cite{mcmullengaussbonnet} using a very different
computation; and by Koziarz and Nguyen~\cite{koziarznguyen} in a
computation which is closer in spirit to ours, since they compute
intersection numbers.

More specifically, in our paper, the orbifold Euler characteristics
are obtained by identifying the quotients as pairs $(X,\Delta)$ where
$X$ is an explicit normal space birational to the quotient of $\P^n$
by a finite group, and $\Delta$ is an explicit $\Q$-divisor in $X$.
We then compute
\begin{equation}\label{eq:hirz}
  \frac{1}{(n+1)^{n-1}}c_1^{orb}(X,\Delta)^n=\frac{(-1)^n}{(n+1)^{n-1}}(K_X+\Delta)^n,
\end{equation}
which is equal to $c_n^{orb}(X,\Delta)$ (and the latter is equal to
the orbifold Euler characteristic). Indeed, by Hirzebruch
proportionality~\cite{hirzprop}, the ratios of Chern numbers for ball
quotients must be the same as those of the compact dual symmetric
space $\P^n$, and we have $c_1(\P^n)=nH$, $c_n(\P^n)=(n+1)H^n$ (where
$H$ denotes the hyperplane class). We will only use the case $n=3$,
where the relevant formula reads $c_1^{orb}(X,\Delta)^3=16c_3(X,\Delta)$.

Note that we do not compute the orbifold Euler characteristic directly,
which can be done by using the stratification of $X$ by strata with
constant isotropy groups (see~\cite{satake}
or~\cite{mcmullengaussbonnet}). Indeed, such a computation would
require a lot of bookkeeping (especially in cases where the relevant
ball quotient is obtained from $\P^3$ by blowing-up and then
contracting, see section~\ref{sec:chl}).

Strictly speaking, the above description is only valid for compact
ball quotients. In terms of the notation in~\cite{chl}, cocompactness
corresponds to the fact that $\kappa_L\neq 1$ for every irreducible
mirror intersection $L$ in the arrangement. On the other hand, the
formulas we use for cocompact lattices remain valid for non-cocompact
ones, since the volume of the complex hyperbolic structures
constructed by Couwenberg, Heckman and Looijenga depend continuously
(in fact even analytically) on the deformation parameter (see
Theorem~3.7 in~\cite{chl}). 

Our computations depend on detailed properties of the
combinatorics of the hyperplane arrangements given by the mirrors in
4-dimensional Shephard-Todd groups. We list these combinatorial
properties in section~\ref{sec:comb} in the form of tables, since we
could not find all of it in the literature (the data can be gathered
fairly easily using modern computer technology).

For concreteness, we also give explicit presentations for the
3-dimensional CHL lattices in terms of generators and relations, see
section~\ref{sec:presentations}. The fact that one can work out
explicit presentations was already mentioned by Couwenberg, Heckman
and Looijenga (see Theorem~7.1 in~\cite{chl}). This depends on the
knowledge of the presentations for braid groups that were worked out
by Brou\'e, Malle, Rouquier~\cite{brmaro}, Bessis and
Michel~\cite{bessismichel}, and fully justified thanks to later work
by Bessis~\cite{bessisannals}.

We hope that our paper provides useful insight into the beautiful
paper by Couwenberg, Heckman and Looijenga.

\noindent \textbf{Acknowledgements:} It is a pleasure to thank
St\'ephane Druel for many useful conversations related to this
paper. I also thank John Parker for his suggestion to use
algebro-geometric methods to compute volumes of CHL lattices. I am
also very grateful to the referee, who made many suggestions to
improve the paper.

\section{Finite unitary groups generated by complex reflections} \label{sec:st}

In this section, we briefly recall the Shephard-Todd classification of
finite unitary groups generated by complex reflections,
see~\cite{shephardtodd} (see also~\cite{cohen} or~\cite{brmaro}).

\subsection{Complex reflections}

Recall that a complex reflection in $V=\C^n$ is a diagonalizable
linear transformation whose eigenvalues are $1,\dots,1,\zeta$ for some
complex number $\zeta\neq 1$ with $|\zeta|=1$ (the eigenvalue 1 has
multiplicity $n-1$). A group is called a \emph{complex reflection
  group} if it is generated by complex reflections, and it is called
unitary if it preserves a Hermitian form on $V$.

Complex reflections preserving a Hermitian inner product $\langle
v,w\rangle=w^*Hv$ can be written as $R_{v,\zeta}$ where
$$
  R_{v,\zeta}(x)=x+(\zeta-1)\frac{\langle x,v\rangle}{\langle v,v\rangle}v 
$$ 
for some nonzero vector $v\in V$. The fixed point set of $R_{v,\zeta}$
in $V$ then consists of the orthogonal complement $v^\perp$ with
respect to $H$, and it is called the \emph{mirror} of $R_{v,\zeta}$.
The number $\zeta$ is called the \emph{multiplier} of $R_{v,\zeta}$
and the argument of $\zeta$ is called the \emph{angle} of
$R_{v,\zeta}$.

We will often assume that the multiplier $\zeta$ is a root of unity
(which is needed if we are to consider only finite groups), and even
that $\zeta=e^{2\pi i/p}$ for some natural number $p\geqslant 2$
(which can be assumed by replacing the reflection by a suitable
power).

\subsection{Braid relations} \label{sec:braiding}

Let $G$ be a group, and let $a,b\in G$. We say that $a,b$ satisfy a
braid relation of length $k$ if 
\begin{equation}\label{eq:braidk}
  (ab)^{k/2}=(ba)^{k/2}.
\end{equation}
When $k$ is odd, $(ab)^{k/2}$ stands for a product of the form
$aba\cdots ba$ with $k$ factors (and similarly for $(ba)^{k/2}$). We
write $\br_k(a,b)$ for the relation of equation~\eqref{eq:braidk}.

Of course, $\br_1(a,b)$ is equivalent to $a=b$, and $\br_2(a,b)$ means
that $a$ and $b$ commute. The relation $\br_3(a,b)$, i.e. $aba=bab$ is
often called the standard braid relation. If $\br_k(a,b)$ holds, but
$\br_j(a,b)$ does not hold for any $j<k$, we write $\br(a,b)=k$.

\subsection{Coxeter diagrams}

It is customary to describe complex reflection groups (with a finite
generating set of reflections) by a Coxeter diagram, which is a
labelled graph. The set of nodes in the graph is given by a
generating set of complex reflections, and each node consists of a
circled integer, corresponding to the order of the corresponding
complex reflection (more precisely, a circled $p$ stands for a complex
reflection of angle $\frac{2\pi}{p}$).

The nodes in the Coxeter diagram, corresponding to reflections $a$ and
$b,$ are joined by an edge labelled by a positive integer $k$ if the
braid relation $\br_k(a,b)$ holds (see~\ref{sec:braiding}). Moreover,
by convention:
\begin{itemize}
\item when $\br_2(a,b)$, the edge is not drawn,
\item when $\br_3(a,b)$, the label 3 is omitted and the corresponding edge is drawn without any label,
\item when $\br_4(a,b)$, the label 4 is omitted and the corresponding edge is drawn as a double edge.  
\end{itemize}

Beware that a given complex reflection group can be represented by
several Coxeter diagrams (since there can be several non-conjugate
generating sets of reflections), and in general a Coxeter diagram need
not represent a unique group (even up to conjugation in $GL(V)$).

\subsection{The Shephard-Todd classification}

Let $G$ be a group acting irreducibly on $V=\C^n$. If $G$ has an
invariant Hermitian form on $V$, then that form is unique. If we
assume further that $G$ is finite, then any invariant Hermitian form
must be definite, so we can think of $G$ as a subgroup of $U(n)$.

Finite subgroups of $U(n)$ generated by complex reflections were
classified by Shephard-Todd in~\cite{shephardtodd}. It is enough to
classify irreducible groups, which come in three infinite families
(symmetric groups, imprimitive groups $G(m,p,n)$ and groups generated
by a single root of unity), together with a finite list of groups.

The infinite families occur in the Shephard-Todd list as $G_1,G_2$ and
$G_3$, and the finite list contains groups $G_4$ through $G_{22}$ (in
dimension 2), $G_{23}$ through $G_{27}$ (in dimension 3), $G_{28}$
through $G_{32}$ (in dimension 4), $G_{33}$ (in dimension 5), $G_{34}$
and $G_{35}$ (in dimension 6), $G_{36}$ (in dimension 7), and $G_{37}$
(in dimension 8). The list is given on p.301 of~\cite{shephardtodd}.

\subsection{Presentations for Shephard-Todd groups and for the associated braid groups}

\begin{figure}
  \includegraphics[width=0.3\textwidth]{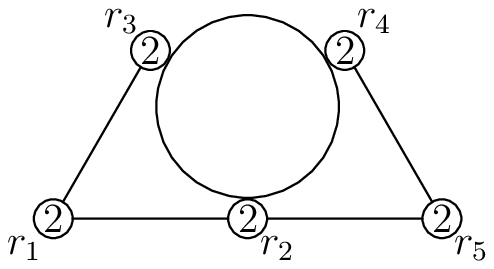}
  \caption{The Coxeter diagram for $G_{31}$}\label{fig:cox-g31}
\end{figure}

Presentations for these groups in terms of generators and relations
are listed in section~11 of~\cite{shephardtodd}. It can be useful to
have reflection presentations, i.e. presentations such that the
generators correspond to complex reflections in the group. Coxeter
diagrams for reflection presentations can be found in convenient form
in pp.185--188 of~\cite{brmaro}; some diagrams have extra decorations,
since the braid relations between generators do not always
suffice. For instance, the diagram for the group $G_{31}$ is the one
given in Figure~\ref{fig:cox-g31}.
This diagram gives a presentation of the form
\begin{eqnarray*}
\langle \quad
r_1,r_2,r_3,r_4,r_5 & | & r_1^2, r_1r_5r_4=r_5r_4r_1=r_4r_1r_5\\  
                    &   &  \br_3(r_1,r_2), \br_3(r_2,r_5), \br_3(r_5,r_3), \br_3(r_3,r_4), [r_1,r_3], [r_2,r_3], [r_2,r_4]
\quad \rangle 
\end{eqnarray*}
where all braid and commutation relations are dictated by the general Coxeter
description, and the circle joining the nodes for $r_1$, $r_5$ and
$r_4$ stands for the relations $r_1r_5r_4=r_5r_4r_1=r_4r_1r_5$.

A more subtle question concerns the presentations of the corresponding
braid groups, which we now define. Given an irreducible finite unitary
group $G$ generated by complex reflections in $V$, we denote by $V^0$
the complement of the union of the mirrors of all complex reflections
in $G$. It is a well known fact that $G$ acts without fixed
points on $V^0$ (this is due to Steinberg~\cite{steinberg}), and the
fundamental group of the quotient $\pi_1(V^0/G)$ is called the braid
group associated to $G$.

For $G=S_{n}$ (acting on $V=\C^{n-1}$ seen as the hyperplane
$\sum z_j=0$ in $\C^{n}$), $\pi_1(V^0/G)$ is simply the usual braid
group $\Br_n$ on $n$ strands (see~\cite{foxneuwirth}). For more
complicated Shephard-Todd groups $G$, presentations were given
in~\cite{brmaro},~\cite{bessismichel}; some of their presentations were
conjectural at the time, but the conjectural statements were later justified
by Bessis~\cite{bessisannals}.
The general rule is that the Coxeter diagrams given in~\cite{brmaro}
give presentations of the corresponding braid groups by removing the
relations expressing the order of the generators. For example, the
braid group associated to $G_{31}$ has the presentation
\begin{eqnarray*}
\langle \quad
r_1,r_2,r_3,r_4,r_5 & | & r_1r_5r_4=r_5r_4r_1=r_4r_1r_5\\  
                    &   &  \br_3(r_1,r_2), \br_3(r_2,r_5), \br_3(r_5,r_3), \br_3(r_3,r_4), [r_1,r_3], [r_2,r_3], [r_2,r_4]
\quad \rangle 
\end{eqnarray*}
We will describe the corresponding group by a diagram whose nodes are
simply bullets without any label giving the order, see the diagrams at
the top of Tables~\ref{fig:comb-A4} through~\ref{fig:comb-g31}
(pp.~\pageref{fig:comb-A4}--\pageref{fig:comb-g31}).

In this paper, we only consider 4-dimensional Shephard-Todd groups
(whose projectivization acts on $\P^3_{\C}$), so we consider
$S_5=W(A_4)$, $G(m,p,3)$ and $G_{28}$ through $G_{32}$.
In fact we will restrict the list even further, because it turns out
some groups will give rise to the same lattices via the
Couwenberg-Heckman-Looijenga construction.

\section{The Couwenberg-Heckman-Looijenga lattices} \label{sec:chl}

We briefly recall some of the results in~\cite{chl}. 
Let $V$ be a finite-dimensional complex vector space. For a subgroup
$G\subset GL(V)$, we denote by $\P G$ the image of $G$ under the
natural map $GL(V)\rightarrow PGL(V)$. Given a
complex linear subspace $\{0\}\subsetneq L\subset V$, we denote by $\P L$
its image in the complex projective space $\P V$.

Now let $G$ be an irreducible finite unitary complex reflection group
acting on the complex vector space $V$, and let $H_i$, $i\in I$ denote
the mirrors of reflections in $G$ (recall that a complex reflection is
a nontrivial unitary transformation which is the identity on a linear
hyperplane, called its mirror). We refer to linear subspaces of the
form $\cap_{j\in J} H_j$ for some $J\subset I$ simply as \emph{mirror
  intersections}.  We denote by $V^0$ the complement of the union of
the mirrors, $V^0=V\backslash \cup_{i\in I} H_i$.

The results in~\cite{chl} produce a family of affine structures on
$V^0$, indexed by $G$-invariant functions $\kappa:I\rightarrow
]0,+\infty[$ such that the holonomy around each mirror $H_j$, $j\in I$
    is given by a complex reflection with multiplier $e^{2\pi
      i\kappa(j)}$. We sometimes denote $\kappa(j)$ by $\kappa_j$ or
    $\kappa_{H}$ when $H=H_j$ for some $j\in I$ (this should cause no
    ambiguity, since we will of course assume $H_i\neq H_j$ when
    $i\neq j$).

For each $\kappa$, up to scaling, there is a unique Hermitian form
which is invariant under the holonomy group.  In what follows, we
assume that the weight assignment $\kappa$ is hyperbolic, in the
sense that the invariant Hermitian form has signature $(n,1)$, where
$\dim V=n+1$.  We denote by $\widetilde{\Gamma}_\kappa$ the holonomy
group, and by $\Gamma_\kappa=\P\widetilde{\Gamma}_\kappa$ its
projectivization.

Couwenberg, Heckman and Looijenga formulate a fairly simple sufficient
condition to ensure that $\Gamma_\kappa=\P\widetilde{\Gamma}_\kappa$
is actually a lattice in $PU(n,1)$, which they refer to as the Schwarz
condition. 

We briefly recall that condition, which is about \emph{irreducible}
mirror intersections in the arrangement (see p.~88 in~\cite{chl} for
definitions). Given a mirror intersection $L=\cap_{j\in J} H_j$, the
set of mirrors containing $L$ induces an arrangement $\mathcal{H}_L$
on $V/L$. We call the mirror intersection $L$ irreducible if
$\mathcal{H}_L$ is irreducible in the sense that it cannot be written
as the product of two lower-dimensional arrangements. Concretely, a
mirror intersection $L$ of codimension $N$ is irreducible if and only
if $N=1$ or there exist $N+1$ mirrors containing $L$ such that $L$ is
the intersection of any $N$ of them.

Given an (irreducible) mirror intersection $L$, we define a real
number $\kappa_L$ as follows. Denote by $G_L$ the fixed point
stabilizer of $L$ in $G$, which is known to be generated by the
reflections in $G$ whose mirror contains $L$ (this follows from
Steinberg's theorem~\cite{steinberg}). Note however that the
stabilizer of $L$ need not be a complex reflection group.

Now define
$$
   \kappa_L=\frac{\sum_{H_j\supset L} \kappa_j}{{\rm codim} L}.
$$ 
With such notation, the Schwarz condition is the requirement that
   for each irreducible mirror intersection $L$ such that
   $\kappa_L>1$,
\begin{equation}\label{eq:schwarz0}
   \frac{|Z(G_L)|}{\kappa_L-1}\in\N,
\end{equation}
where $Z(G_L)$ denotes the center of $G_L$. 
\begin{rk}
  This condition is a generalization of the Mostow $\Sigma$-INT
  condition~\cite{mostowihes}. The analogue of the INT condition
  in~\cite{delignemostow} would be the requirement that
  $(\kappa_L-1)^{-1}\in \N$. For more on this, see Example~4.3, p.~131
  of~\cite{chl}.
\end{rk}

Applied to the case where $L$ is a single mirror $H=H_i$, fixed by a
reflection of maximal order $o$ in $G$, the Schwarz condition says
that
\begin{equation}\label{eq:defint}
  \kappa_H=1-\frac{o}{p_H}
\end{equation}
for some integer $p_H$. In fact, it is enough to consider
Shephard-Todd generated by reflections of order 2 (the other ones
would not produce any more lattices in $PU(n,1)$), in which case
condition~\eqref{eq:defint} reads $\kappa_H=1-\frac{2}{p_H}$.

\begin{rk}The holonomy group acts irreducibly on $V=\C^{n+1}$, and preserves a
unique Hermitian form. The signature of the Hermitian form can be read
off the number $\kappa_{\{0\}}$, namely the form is definite when
$\kappa_{\{0\}}<1$, degenerate for $\kappa_{\{0\}}=1$, and hyperbolic
for $\kappa_{\{0\}}>1$.
\end{rk}

Couwenberg-Heckman-Looijenga lattices are indexed by a Shephard-Todd
group $G$ and a list of integers, one for each $G$-orbit of mirror of
a complex reflection in $G$ (these integers are the ones in
equation~\eqref{eq:defint}, one from each $G$-orbit of mirror). For
most groups $G$, there is a single orbit of mirrors, and we denote the
corresponding integer by $p$ , and the projectivized holonomy group by
$\CHL(G,p)$.  In some cases, there are two $G$-orbits of mirrors, in
which case we denote the two integers by $p_1,p_2$, and the group by
$\CHL(G,p_1,p_2)$ (we use a natural order in the $G$-orbits of mirrors
corresponding to a numbering of the generators, see
section~\ref{sec:comb}). To get a uniform notation for both cases, we
will denote the group by $\CHL(G,{\bf p})$, where ${\bf p}$ stands for
either $p$ or $p_1,p_2$.

An important result in~\cite{chl} is the following.
\begin{thm}
  Suppose that $\kappa$ a hyperbolic $G$-invariant function such that
  the Schwarz condition is satisfied, and denote by ${\bf p}$ the
  integers coming from equation~\eqref{eq:defint}. Then the group
  $\CHL(G,{\bf p})$ is a lattice in $PU(n,1)$.
\end{thm}

The arithmeticity of the corresponding lattices was studied
in~\cite{chl3d_1}. In dimension at least three, the list turns out to
contain only two commensurability classes of non-arithmetic lattices,
both in $PU(3,1)$.

In dimension $n>1$, there are only finitely many choices of ${\bf p}$
such that the Schwarz condition is satisfied, which are listed
in~\cite{chl}, pp.~157--160 (see also~\cite{chl3d_1} and
section~\ref{sec:comb} of this paper). In order to produce the list,
one needs to known some detailed combinatorial properties of the
arrangements, which are listed in section~\ref{sec:comb} of our paper.

The list contains the Deligne-Mostow lattices (for which the Schwarz
condition is equivalent to the generalized Picard integrality
condition) and the Barthel-Hirzebruch-H\"ofer lattices in $PU(2,1)$.

In this paper, we list only 3-dimensional groups (the corresponding
finite unitary groups act on $\C^4$). Moreover, we only consider
$G$-invariant weight assignments $\kappa$ (for most arrangements, any
$\kappa$ for which the more general Couwenberg-Heckman-Looijenga
results apply are actually $G$-invariant, see section~2.6
in~\cite{chl}, in particular Proposition~2.33). In particular, we do
not reproduce the entire Deligne-Mostow list (which contains many non
$G$-invariant assignments).

In order to prove their result, Couwenberg, Heckman and Looijenga
consider the developing map of their geometric structures, which is a
priori only defined on an unramified covering $\widetilde{V}^0$ of
$V^0$ (the holonomy covering, which is the covering whose fundamental
group is the kernel of the holonomy representation). They show that
the developing map extends above a suitable blow-up $\widehat{V}$ of
$V$, namely the one obtained by blowing-up linear subspaces
corresponding to mirror intersections $L$ with $\kappa_L>1$, in order
of increasing dimension (this can be done in a $G$-invariant manner,
since $\kappa$ is assumed $G$-invariant). We denote by $\widehat{X}$
the corresponding blow-up of projective space $X=\P(V)$.

The proof of the Couwenberg-Heckman-Looijenga results require careful
analysis of where the developing map is a local biholomorphism, which
does not usually happen on all the exceptional divisors of the blow-up
$\widehat{X}$. In fact the components above exceptional divisors
$D(L)$ corresponding to a mirror intersection $L$ of (linear)
dimension $k$, get mapped to components of \emph{codimension} $k$
under the developing map (see part~(ii) of Proposition~6.9
in~\cite{chl}).

Since we consider only 3-dimensional ball quotients, we will need to
handle only situations where the $L$'s that get blown-up to obtain
$\widehat{V}$ have dimension 1 or 2 (i.e. these correspond to points
or lines in the projective arrangement in $\P V$).

When blowing a point in $\P V$ that corresponds to a mirror
intersection $L$ which is a line (but we do not blow up any
higher-dimensional mirror intersection that contains it), the
developing map is a local biholomorphism above that point, since the
corresponding exceptional divisor gets mapped to a divisor (see
Proposition~6.9 in~\cite{chl}).

A slightly more complicated situation occurs when, among the mirror
intersections, there is a 2-plane $L$ such that $\kappa_L>1$ (this
corresponds to a line in $\P(V)$). In that case, every irreducible
1-dimensional mirror intersection $M$ with $M\subset L$ also satisfies
$\kappa_M>1$, see the monotonicity statement in Corollary~2.17
of~\cite{chl} (note that such $M$ correspond to points $\P M$ in
$\P(V)$).  In particular, $\widehat{V}$ is obtained by first
blowing-up all the 1-dimensional mirror intersections contained in
$L$, then blowing-up the strict transform of $L$ (see
Figure~\ref{sfig:bu2} on p.~\pageref{sfig:bu2} where we have blown up
points, and~\ref{sfig:bu3} where we have blown up the strict tranform
of the lines joining these points). We denote by $D(M)$ and $D(L)$ the corresponding
exceptional divisors in $\widehat{X}$.

The developing map then maps $D(M)$ to a divisor, and $D(L)$ to a
variety of codimension 2, i.e. a curve. In particular, in order to
describe the corresponding ball quotient, we will need to contract the
divisors $D(L)$ to curves. It turns out the $D(L)$ we will encounter
(still with $\dim L=2$) are actually copies of $\P^1\times \P^1$, see
p.~150 in~\cite{chl}. There are two ways to contract them, by
collapsing one or the other factor. Of course, since the developing
map does not extend to $\widehat{X}$, we will contract in the
direction opposite to the one that gives $\widehat{X}$ (see
how Figure~\ref{sfig:bu3} collapses to Figure~\ref{sfig:bu4}).

Note once again that the above blow-up can be performed in a
$G$-invariant manner, since our weight assignement is assumed to be
$G$-invariant.

In our volume formulas, we will need a description of the canonical
divisors of $K_{\widehat{X}}$ and $K_Y$. The first remark is that $Y$
is a normal space, so $K_Y$ is defined, and it has
$\mathbb{Q}$-factorial singularities. This follows from Lemma~5.16
in~\cite{kollarmori}, since $Y$ is a quotient of the unit ball by a
lattice (see Proposition~\ref{prop:key}), and lattices have normal
torsion-free subgroups of finite index, so $Y$ is the quotient of a
(quasi-)projective algebraic variety by a finite group.

In all cases we consider, the strict transforms of the lines that
we blow up are pairwise disjoint, so it is enough to consider the case
where we blow up the strict transform of a single line. Consider a
line $L$ in $\P^3$, and denote by $\pi_1$ the blow-up of $n$ distinct
points on $L$ (we assume $n\geq 2$). Denote by $\pi_2$ the blow-up of
the strict transform of $L$ and by
$\pi=\pi_1\circ\pi_2:\widehat{X}\rightarrow \P^3$ the composition.  We
denote by $D_1,\dots,D_n$ the exceptional locus over the points that
were blown-up in $\pi_1$, and by $E$ the exceptional divisor over $L$.

Note that $E$ is isomorphic to $\P^1\times\P^1$,
in particular $Pic(E)\simeq \Z l_1\oplus \Z l_2$, where we assume
$l_1$ projects to $L$ in $\P^3$. We then have
\begin{equation}\label{eq:self}
  E|_E=-l_1-(n-1)l_2,\quad E^3=(-l_1-(n-1)l_2)^2=2n-2.
\end{equation}
We will be interested in the space obtained from $\widehat{X}$ by
contracting $E$ to a $\P^1$, by contracting the factor given by
$l_1$. We denote by $f:\widehat{X}\rightarrow Y$ the corresponding
map.

Note once again that the space $Y$ is singular (unless $n=2$), but it
has $\Q$-factorial singularities, which implies that we can pull-back
the canonical divisor $K_Y$ under $f$, and we have
\begin{equation}\label{eq:discrepancy0}
K_{\widehat{X}}=f^*K_Y+aE
\end{equation}
for some rational number $a$.

Now we take the intersection of both sides of
equation~\eqref{eq:discrepancy0} with $l_1$, and note that
$$
E\cdot l_1=E|_E\cdot l_1=-(n-1).
$$
Using the adjunction formula $K_E=K_{\widehat{X}}|_E+E|_E$, we have
$$
K_{\widehat{X}}\cdot l_1=K_{\widehat{X}}|_E\cdot l_1=(K_E-E|_E)\cdot l_1=n-3,
$$
so we get $a=-(n-3)/(n-1)$ and
\begin{equation}\label{eq:discrepancy}
  K_{\widehat{X}}=f^*K_Y-\frac{n-3}{n-1}E.
\end{equation}

We will also need to study $f^*f_*Z$ for various divisors $Z$. The
following follows from computations similar to the previous one.
\begin{prop}\label{prop:formula}
  Let $Z\subset\widehat{X}$ be the proper transform of a plane in
  $H\subset \P^3$.
  \begin{enumerate}
    \item If $H\cap L$ is one of the points blown-up in $\pi_1$, or if
      $H$ contains $L$, then $f^*f_*Z=Z$.
    \item If $H\cap L$ is a point which is not one of the points blown-up in $\pi_1$, then
      $$
      f^*f_*Z=Z+\frac{1}{n-1}E.
      $$
    \item For every $j$, we have
      \begin{equation*} 
        f^*f_*D_j=D_j+\frac{1}{n-1}E.
      \end{equation*}
  \end{enumerate}
\end{prop}
Parts~(2) and ~(3) follow from the fact that (under the hypothesis on
$H$) $Z|_E$ and $D_j|_E$ are equivalent to $l_2$.

%

The log-canonical divisor corresponding to the CHL ball quotient will
be $K_Y+\Delta$, where $\Delta$ denotes
\begin{equation}\label{eq:Delta}
  \Delta = \sum_{i\in I} \kappa_i f_*\widehat{H}_i + \sum_{\dim
    L=1,\kappa_L>1} \left(2-\kappa_L\right) f_*D(L).
\end{equation}
In formula~\eqref{eq:Delta}, $L$ ranges over all irreducible mirror
intersections, and $\widehat{H}_i$ denotes the strict transform of
$H_i$ under the blow-up $\pi:\widehat{X}\rightarrow X$.

Since our function $\kappa$ is $G$-invariant, the finite
complex reflection group $G$ acts on $\widehat{X}$ and on $Y$. We
denote by $\varphi:Y\rightarrow Y/G$ the quotient map, and by
$\mathcal{D}=\varphi_*\Delta$. The irreducible components of
$\mathcal{D}$ correspond to the $G$-orbits of mirror intersections
$L$ as in the sum~\eqref{eq:Delta}. Moreover, $\varphi$ ramifies to
the order $|Z(G_L)|$ around these components, and the coefficients of
$\mathcal{D}$ have the form $1-1/n_L$, where $n_L$ is the integer that
occurs in the Schwarz condition~\eqref{eq:schwarz0}.


We now formulate a key result of Couwenberg, Heckman and Looijenga as
follows (see Theorem~6.2 of~\cite{chl}).
\begin{prop}\label{prop:key}
  Suppose the weight assignment $\kappa$ satisfies the Schwarz
  condition, and denote by ${\bf p}$ the corresponding integers
  attached to $G$-orbits of mirrors. 
  \begin{enumerate}
    \item If $\kappa_L\neq 1$ for every irreducible mirror
      intersection $L$, then the lattice $\CHL(G,{\bf p})$ is
      cocompact, and the quotient $\B^n/\CHL(G,{\bf p})$ is given as
      an orbifold by the pair $(Y/G,\mathcal{D})$;
    \item Otherwise, the ball quotient $\B^n/\CHL(G,{\bf p})$ has one
      cusp for each $G$-orbit of irreducible mirror intersection $L$
      with $\kappa_L=1$, and it is given as an orbifold by the pair
      $(Y^0/G,\mathcal{D}^0)$, where $Y^0$ is obtained from $Y$ by
      removing the image of the irreducible mirror intersections $L$
      with $\kappa_L=1$.
  \end{enumerate}
\end{prop}
\begin{rk}
  If we require the stronger condition $(\kappa_L-1)^{-1}\in\N$
  instead of the Schwarz condition~\eqref{eq:schwarz0}, then
  $(Y,\Delta)$ a also ball quotient orbifold, that covers
  $(Y/G,\mathcal{D})$. In other words, $\CHL(G,{\bf p})$ then has a
  sublattice of index $|\P G|$, which is the orbifold fundamental
  group of $(Y,\Delta)$.
\end{rk}

As discussed in the introduction, the volume of the quotient can be
computed up to a universal multiplicative constant as the
self-intersection
$$
  c_1(Y/G,\mathcal{D})^3=\frac{c_1(Y,\Delta)^3}{|\P G|}=-\frac{(K_Y+\Delta)^3}{|\P G|}.
$$
We will work out several specific examples of this general
construction in section~\ref{sec:volumes}.

Note that a lot of the above description makes sense when the Schwarz
condition is not satisfied. If the weight assignment $\kappa$ is
hyperbolic, one gets a complex hyperbolic cone manifold structure on
$(Y,\Delta)$, but the coefficients in the divisor
$\mathcal{D}=f_*\Delta$ are no longer of the form $1-1/k$ for $k$ an
integer, so $(Y/G,\mathcal{D})$ is not an orbifold pair.
It follows from Theorem~3.7 in~\cite{chl} that the volume of these
structures depends continuously on the parameters ${\bf p}$ (see
equation~\eqref{eq:defint}), because of the analyticity of the
dependence on ${\bf p}$ of the Hermitian form invariant by the
holonomy group. Indeed, the Riemannian metric can be expressed in
terms of the Hermitian form, see p.~135 in~\cite{mostowbook}, and the
volume form is simply the square-root of the determinant of the
corresponding metric. The volume can then be computed as an integral
on a possibly blown-up projective space (the blow-up does not affect
volume since it is an isomorphism away from a set of measure zero).

In particular, in order to compute the volume of a non-compact ball
quotient for some parameter $p_H$ (i.e. one where $\kappa_L=1$ for
some $L$), one can compute the volume of the structures for
$p_H-\varepsilon$, then let $\varepsilon$ tend to 0. Some of the
methods below require $\Delta$ to be a $\Q$-divisor, so we should
actually take $\varepsilon$ rational. The upshot is that our volume
computations are valid even for non-cocompact lattices, and we will
not need to consider the Couwenberg-Heckman-Looijenga compactification
in those cases.

\section{Relation with Deligne-Mostow groups} \label{sec:dm}

As mentioned in~\cite{chl} (see their section~6.3), the
Couwenberg-Heckman-Looijenga construction applied to reflection groups
of type $A_n$ and $B_n$ give lattices commensurable to the
Deligne-Mostow examples. We give some details of that relationship in
the case of lattices in $PU(3,1)$, in the form of a table (see
Figure~\ref{tab:b4-dm}). The basic point is that each Deligne-Mostow
lattice in $PU(n,1)$ is the image of a representation of a spherical
braid groups on $N=n+3$ strands (which is isomorphic to the
corresponding plane braid group $\Br_N$ modulo its center), and the
representation is determined by the choice of an $N$-tuple
$\mu=(\mu_1,\dots,\mu_N)$ of hypergeometric exponents, and a subgroup
$\Sigma\subset S_N$ that leaves $\mu$ invariant.

The group that gets represented is $\phi^{-1}(\Sigma)$ for some
subgroup $\Sigma\subset S_N$ ($\Sigma$ acts as a symmetry group of the
$N$-tuple of weights for the corresponding hypergeometric functions),
where $\phi:\Br_N\rightarrow S_N$ corresponds to remembering only the
permutation effected by the braid. The corresponding hypergeometric
group is denoted by $\Gamma_{\mu,\Sigma}$ (see~\cite{mostowihes}
or~\cite{mostowsurvey}).

For simplicity, when describing Deligne-Mostow groups, we will take
$\Sigma$ to be the full symmetry group of the $N$-tuple
$\mu=(\mu_1,\dots,\mu_N)$ of weights, but the corresponding CHL
subgroups will be obtained by taking $\Gamma_{\mu,\Sigma_0}$ for some
subgroups $\Sigma_0\subset\Sigma$.  For the reader's convenience, the
explicit commensurabilities are summarized in Table~\ref{tab:b4-dm},
on p.~\pageref{tab:b4-dm}. We briefly explain how to obtain this table in
sections~\ref{sec:5fold} and~\ref{sec:4fold}.

\subsection{Hypergeometric monodromy groups with 5 equal exponents}\label{sec:5fold}

Suppose first that $\mu=(\mu_1,\dots,\mu_6) \in ]0,1[^6$ has 5 equal
values, say $\mu_1=\mu_2=\dots=\mu_5=\alpha$. We assume moreover that
$\mu$ satisfies condition $\Sigma$-INT for $\Sigma=S_5$, in particular
$1-2\alpha=2/p$ for some integer $p$. The sixth exponent $\mu_6$ is
determined by $\alpha$, since $\sum_{j=1}^6 \mu_j=2$, so
$\mu_6=2-5\alpha=-\frac{1}{2}+\frac{5}{p}$. In particular, we must
have $3<p<10$, since we want $0<\mu_6<1$.

Condition $\Sigma$-INT requires that, if $\alpha+\mu_6<1$, then
$(1-\alpha-\mu_6)^{-1}=1/k$ for some $k\in\N\cup\{\infty\}$. This
implies $p=4,5,6$ or $8$ (see the top half of
Table~\ref{tab:b4-dm}).

\begin{figure}
  \includegraphics[width=0.3\textwidth]{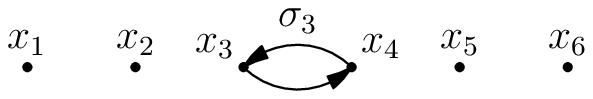}
  \caption{A half-twist, standard generator $\sigma_3$ of the planar
    braid group on 6 strands.} \label{fig:braids}
\end{figure}

Consider the standard generators of the braid group, i.e half-twists
$\sigma_j$ between $x_j$ and $x_{j+1}$ ($j=1,\dots,5$), see
Figure~\ref{fig:braids}. These satisfy the well known standard braid
relation
$\sigma_j\sigma_{j+1}\sigma_j=\sigma_{j+1}\sigma_{j}\sigma_{j+1}$, and
$\sigma_j$ commutes with $\sigma_k$ for $|j-k|\geqslant 2$.

In fact these relations give a presentation for the braid group
$\Br_6$ corresponding to the sixtuples of points in $\C$ (this was
proved by Artin, see~\cite{artin}). From this, one can deduce a presentation
of the corresponding spherical braid group, i.e. the one corresponding
to sixtuples of points in $P^1_\C$. One verifies in particular that
(the images in the spherical braid group of) the first four generators
$\sigma_1,\dots,\sigma_4$ suffice to generate the group, see the
discussion in~\cite{mostowsurvey}.

It is well known that the monodromy of the $\sigma_j$ are complex
reflections, with non-trivial eigenvalue $e^{2\pi i/p}$
(see~\cite{delignemostow} or~\cite{thurstonshapes}), so the
hypergeometric monodromy group is a homomorphic image of the braid
group $\Br_5$, such that the four standard generators are mapped to
reflections of angle $2\pi/p$.  We then get the following.
\begin{prop}
  Let $p=4,5,6$ or $8$,
  $\mu=(\frac{1}{2}-\frac{1}{p},\frac{1}{2}-\frac{1}{p},\frac{1}{2}-\frac{1}{p},\frac{1}{2}-\frac{1}{p},\frac{1}{2}-\frac{1}{p},\frac{5}{p}-\frac{1}{2})$. Then
  the group $\Gamma_{\mu,S_5}$ is conjugate in $PU(3,1)$ to
  $\CHL(A_4,p)$.
\end{prop}

In the case $p=6$, the sixtuple $\mu=(1,1,1,1,1,1)/3$ has a larger
symmetry group, namely $S_6$ instead of $S_5$. If we write
$\Sigma=S_6$ and $\Sigma_0=S_5$, then $\Gamma_{\mu,\Sigma_0}$ has
index 6 in $\Gamma_{\mu,\Sigma}$ (because $\Sigma_0$ has index 6 in
$\Sigma$, and $\mu$ satisfies condition INT, see~\cite{mostowihes}).

\subsection{Hypergeometric monodromy groups with 4 equal exponents}\label{sec:4fold}

\begin{figure}
  \includegraphics[width=0.3\textwidth]{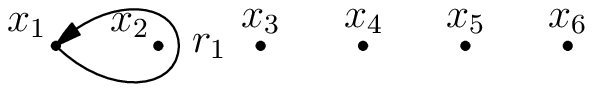}
  \caption{A full twist between $x_1$ and $x_2$.} \label{fig:braids2}
\end{figure}

Suppose now the sixtuple of exponents $\mu$ has 4 equal exponents, say
$\mu_2=\mu_3=\mu_4=\mu_5$, and satisfies the $\Sigma$-INT condition
with respect to $\Sigma\simeq S_4$. As in the previous section, the
monodromy group is generated by $r_1=\sigma_1^2$, $r_2=\sigma_2$,
$r_3=\sigma_3$ and $r_4=\sigma_4$. The loop $r_1$ is called a
full twist between $x_1$ and $x_2$, see Figure~\ref{fig:braids2}.

Now it is easy to see that if two group elements $a,b$ satisfy the
braid relation $aba=bab$, then $c=a^2$ and $d=b$ satisfy a higher
braid relation, namely $(cd)^2=(dc)^2$ (see section~\ref{sec:braiding}
for the definition of braiding). The subgroup of $\Br_6$ generated by
$r_1,\dots,r_4$ is then a braid group of type $B_4$, see the Coxeter
diagram of Figure~\ref{fig:comb-B4} (p.~\pageref{fig:comb-B4}).

This shows that the corresponding Deligne-Mostow group
$\Gamma_{\mu,S_4}$ is a homomorphic image of the braid group of type
$B_4$. Moreover, it generated by complex reflections of understood
angles, namely $\rho(r_1)=\rho(\sigma_1^2)$ rotates by
$2\pi(1-\mu_1-\mu_2)$, whereas for $j=2,3,4$,
$\rho(r_j)$ rotates by $2\pi(\frac{1}{2}-\mu_2)$.  The
$\Sigma$-INT condition says that these two angles can be written as
$2\pi/p_1$ and $2\pi/p_2$ respectively, where $p_1,p_2$ are integers.
We then have:
\begin{prop}
  With the above notation, the group $\Gamma_{\mu,S_4}$ is conjugate
  in $PU(3,1)$ to $\CHL(B_4,p_1,p_2)$.
\end{prop}
The list of $\mu$ and the corresponding pairs $(p_1,p_2)$ is given in
Table~\ref{tab:b4-dm}.
For example, for the Deligne-Mostow group $\Gamma_{\mu,\Sigma}$ with
$\mu=(5,3,3,3,3,7)/12$, $\Sigma=S_4$, we get $p_1=3$ and $p_2=4$, or
in other words $\Gamma_{\mu,S_4}$ is conjugate to $\CHL(B_4,3,4)$. Of
course, we can switch the exponents $\mu_1$ and $\mu_6$, which gives
another description of $\Gamma_{\mu,S_4}$ as $\CHL(B_4,6,4)$.

\section{Computation of the volumes}\label{sec:volumes}

As in section~\ref{sec:chl}, we denote by $\pi:\widehat{X}\rightarrow X$
the blow-up, and by $f:\widehat{X}\rightarrow Y$ the corresponding
contraction, see diagram~\eqref{eq:diag}.
\begin{equation}\label{eq:diag}
  \begin{tikzcd}
    & M,D,E\subset \widehat{X}\arrow{dl}{\pi}\arrow{dr}{f} & \\
\pi_*M\subset X &                                               &  f_*M,f_*D\subset Y
  \end{tikzcd}
\end{equation}
We denote by $D$ (resp. $E$) the exceptional locus corresponding to
blowing up irreducible mirror intersections that are points
(resp. lines), and $M$ is the proper transform in $\widehat{X}$ of the
arrangement in $\P^3$.

If there is only one $G$-orbit of mirrors and one orbit of
1-dimensional mirror intersection, then the divisor $\Delta$ reads
$$
  \left(1-\frac{2}{p}\right)f_*M+\left(1-\frac{1}{m}\right)f_*D,
$$ 
where $p$ is the order of the complex reflections for the holonomy
around a hyperplane (more precisely the non-trivial eigenvalue is
$e^{2\pi i/p}$), and $m$ is a rational number computed from $\kappa_L$
as in section~\ref{sec:chl}.

If the arrangement has more than one $G$-orbit of mirrors, we write
$M=\Sigma_j M_j$, and replace $(1-\frac{2}{p})f_*M$ by a sum
$\sum_j(1-\frac{2}{p_j})M_j$ (it turns out in all CHL examples, there
are at most two mirror orbits, i.e. the sum has at most two terms).  A
similar remark is of course in order for $f_*D$ and $E$, since in
general we may have to blow up several $G$-orbits of mirror
intersections.

In the next few sections, we go through the computations of volumes
for CHL lattices associated to the primitive 4-dimensional
Shephard-Todd groups (the corresponding projective space has dimension
3).
The results of the volume computations are given in
Tables~\ref{tab:b4-dm} and~\ref{fig:rough_invariants} on
pp.~\pageref{tab:b4-dm}--\pageref{fig:rough_invariants}.

In section~\ref{sec:a4}, we treat the case of lattices obtained by the
Couwenberg-Heckman-Looijenga construction from the Weyl group of $A_4$
in detail. The covolumes of the corresponding lattices are known
(see~\cite{mcmullengaussbonnet}), since they are commensurable to
specific Deligne-Mostow lattices, see section~\ref{sec:dm}.  The
corresponding arrangement can be visualized, which should help the
reader follow the computations (first in this simple case, then in
more complicated ones). Indeed $W(A_4)$ is a Coxeter group, the
corresponding arrangement is real, and it contains only 10
hyperplanes, so we can draw a picture, see Figure~\ref{sfig:bu1}.

In subsequent sections, we will treat the groups derived from
$G_{28}$, $G_{29}$, $G_{30}$, $G_{31}$ and $W(B_4)$, where we cannot
draw pictures. The combinatorial properties of these arrangements are
listed in Figures~\ref{fig:comb-A4} through~\ref{fig:comb-g31}
(pp.~\pageref{fig:comb-A4}-\pageref{fig:comb-g31}).

There is an extra (primitive, irreducible) 4-dimensional Shephard-Todd
group, namely $G_{32}$, but just as in~\cite{chl} we omit it from the
list, since $\P(G_{32})=\CHL(W(A_4),3)$, so $G_{32}$ would produce the
same list of complex hyperbolic lattices as $W(A_4)$.

\subsection{The groups derived from the $A_4$ arrangement} \label{sec:a4}

\begin{figure}
  \begin{subfigure}{0.2\textwidth}
    \includegraphics[width=\textwidth]{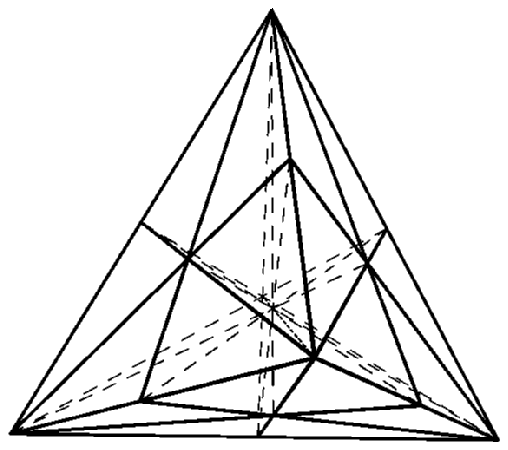}
    \caption{$p=4$}\label{sfig:bu1}
  \end{subfigure}\hfill
  \begin{subfigure}{0.2\textwidth}
    \includegraphics[width=\textwidth]{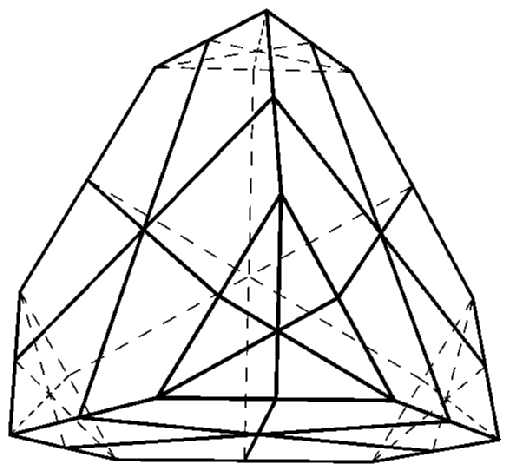}
    \caption{$p=5,6$}\label{sfig:bu2}
  \end{subfigure}\hfill
  \begin{subfigure}{0.2\textwidth}
    \includegraphics[width=\textwidth]{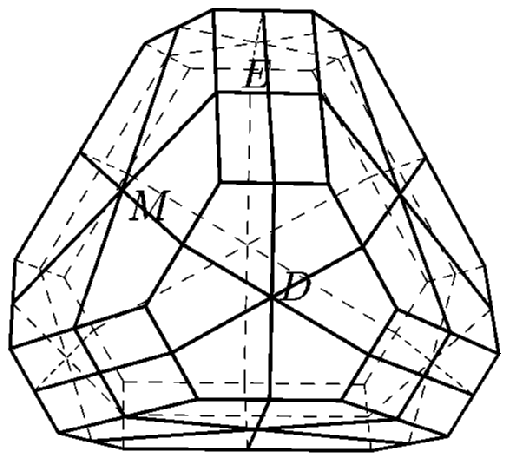}
    \caption{}\label{sfig:bu3}
  \end{subfigure}\hfill
  \begin{subfigure}{0.2\textwidth}
    \includegraphics[width=\textwidth]{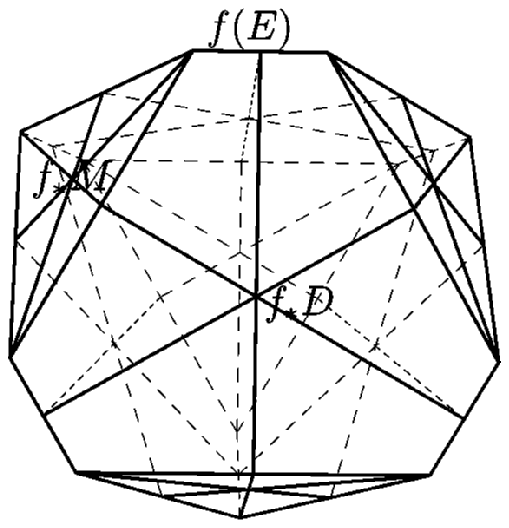}
    \caption{$p=8$}\label{sfig:bu4}
  \end{subfigure}
  \caption{A schematic picture of the process of blowing-up points,
    then lines, and finally contracting $E$ in the other
    direction.}\label{fig:blowup}
\end{figure}

A schematic picture of the projectivization of the $A_4$ arrangement
appears in Figure~\ref{sfig:bu1} (we draw the picture in an affine
chart $\R^3\subset P^3_{\R}\subset P^3_{\C}$). It can be thought of as
the barycentric subdivision of a tetrahedron, but there more symmetry
than the usual euclidean symmetry of the tetrahedron, since the
vertices of the tetrahedron can actually be mapped to the barycenter
via an element of $W(A_4)$.

The group $W(A_4)$ acts transitively on the set of mirrors, so only
one parameter $p$ is allowed, and the weight function $\kappa$ is
constant equal to $1-2/p$ (see equation~\eqref{eq:defint}).

The Schwarz condition holds precisely for $p=2,3,4,5,6$ and $8$ (see
Figure~\ref{fig:comb-A4}, p.~\pageref{fig:comb-A4} for the values of
$\kappa_L$ for various strata $L$). For $p=2$ we recover $W(A_4)$
itself, and for $p=3$ we obtain the Shephard-Todd group $G_{32}$. In
particular, the CHL construction applied to the group $G_{32}$ would
produce the same list of lattices as the one for $W(A_4)$ (more
precisely, the $G_{32}$ arrangement with constant weight function
$1-3/q$ gives the same complex hyperbolic structures as the $A_4$
arrangement with constant function $1-2/q$).

For $p=4,5,6$ or $8$, we get lattices in $PU(3,1)$, which we denote by
$\CHL(A_4,p)$. The volume computation depends on detailed
combinatorial properties of the weighted arrangement, and in
particular the volume formulas depend on $p$ (see
sections~\ref{sec:a4-noblowup} through~\ref{sec:a4-lines}).
What we need from the combinatorics is listed in
Figure~\ref{fig:comb-A4} (on p.~\pageref{fig:comb-A4}).

\subsubsection{The group $\CHL(A_4,4)$} \label{sec:a4-noblowup}

The computation is very easy in this case, since there is no
irreducible mirror intersection $L$ with $\kappa_L>1$. This means that
we do not need any blow-up, i.e. $X=\widehat{X}=Y$, and the orbifold
locus is supported by the hyperplane arrangement.

The arrangement has 10 hyperplanes, so the log-canonical divisor is
numerically equivalent to $(-4+10(1-\frac{2}{4}))H$, where $H$ denotes
the class of a hyperplane, so $(K_X+D)^3=1$.  We have $|\P
G|=|G|=120$, and 3-dimensional ball quotients satisfy
$c_1(X)^3=16c_3(X)$, so the Euler characteristic is given by
$$
-\frac{1}{120\cdot 16}=-\frac{1}{1920}.
$$ 
This is the orbifold Euler characteristic of the Deligne-Mostow
lattice for hypergeometric exponents $\mu=(1,1,1,1,1,3)/4$, see Table~3
in~\cite{mcmullengaussbonnet}.

\subsubsection{The groups $\CHL(A_4,5)$ and $\CHL(A_4,6)$} \label{sec:a4-points}

In the case $p=5$, we have $\kappa_{L_{12}}=1-\frac{1}{10}$ and
$\kappa_{L_{123}}=1+\frac{1}{5}>1$, so we need to blow up the five
points in the $G$ orbit of $\P L_{123}$. Let
$\pi:\widehat{X}\rightarrow X$ denote that blow up.
A schematic picture of the blow up is given in Figure~\ref{sfig:bu2}
with some inaccuracy in the representation, because the barycenter of
the tetrahedron, which is in the same $W(A_4)$ orbit as the vertices,
should be blown-up as well (but this would be too cumbersome to draw).

We denote by $M=M_1+\dots+M_{10}$ the proper transform in
$\widehat{X}$ of the arrangement. Since there are 6 mirrors through
(every element in the $G$-orbit of) $L_{123}$, we have
$$
\pi^*\pi_*M = M + 6 D,\quad K_{\widehat{X}}=\pi^*K_X+2D.
$$ 
In the last formula, the factor of $2$ comes from the codimension
minus one for the locus blown-up in $\P^3$.  We then compute, for
$\nu,\delta\in\Q$,
\begin{eqnarray*}
&\left(K_{\widehat{X}}+\nu M+\delta D\right)^3\\
&=\left( \pi^*K_X+\nu M+(2+\delta)D\right)^3\\
&=\left( \pi^*(K_X+\nu\pi_*M) +\alpha D\right)^3\\
&=\lambda^3 + 5\alpha^3,
\end{eqnarray*}
where
\begin{eqnarray*}
&\lambda=-4+10\nu\\
&\alpha=2+\delta-6\nu.
\end{eqnarray*}
Specializing to $\nu=1-\frac{2}{p}$,
$\delta=2-\kappa_{L_{123}}=1-\frac{1}{5}$ (see the tables in
section~\ref{sec:comb}), we get $(c_1^{orb})^3=136/25$, hence 
$$
\chi^{orb}=c_3^{orb}=-\frac{136}{25\cdot 120\cdot 16}=-\frac{17}{6000}.
$$ 
This agrees with the Euler characteristic of the Deligne-Mostow
lattice for with $\mu=(3,3,3,3,3,5)/10$.

The same formula also works for $p=6$, where we take $\nu=1-\frac{2}{6}$,
$\delta=2-\kappa_{L_{123}}=1-\frac{1}{3}$.  In that case, we get
$$
\chi^{orb}=-\frac{1}{270}.
$$ 
This is coherent with the formula in~\cite{mcmullengaussbonnet},
note that this lattice has index $6=6!/5!$ in the corresponding
Deligne-Mostow lattice, i.e. the one with $\mu=(1,1,1,1,1,1)/3$, and
$(-1/270)/6=-1/1620$ (see section~\ref{sec:dm} for the relationship
with Deligne-Mostow lattices).

\subsubsection{The group $\CHL(A_4,8)$} \label{sec:a4-lines}

In this section, we treat the group derived from $A_4$ with $p=8$,
which corresponds to the Deligne-Mostow group for
$\mu=(1,3,3,3,3,3)/8$. Recall that the combinatorics of the
arrangement are given in Figure~\ref{fig:comb-A4}
(p.~\pageref{fig:comb-A4}), see also Figure~\ref{sfig:bu1}.

The irreducible mirror intersections $L$ with $\kappa_L>1$ consist of
the 5 lines in the $G$-orbit of $L_{123}$ (in the schematic picture,
these correspond to the vertices of the tetrahedron as well as its
barycenter), and the 10 two-planes in the $G$-orbit of $L_{12}$ (these
correspond to the 6 edges of the tetrahedron, together with its 4
lines joining a vertex to the barycenter).  We write $X=\P^3$,
$\pi_1:\widetilde{X}\rightarrow X$ for the blow-up of the five points,
$\pi_2:\widehat{X}\rightarrow\widetilde{X}$ for the blow-up of the
strict transform of the 10 lines through these 5 points, and finally
$\pi=\pi_1\circ\pi_2$.  We use the notation from
section~\ref{sec:volumes}, so that $D$ (resp. $E$) denotes the
exceptional divisor in $\widehat{X}$ above points (resp. lines) in
$X=\P^3$. The space $\widehat{X}$ is depicted in Figure~\ref{sfig:bu3}
(except that we omit drawing the blown-up barycenter and the
exceptional divisors above lines incident to the
barycenter).

The components of $E$ are copies of $\P^1\times\P^1$, and the space
$Y$ is obtained from $\widehat{X}$ by contracting the fibers of these
copies of $\P^1\times\P^1$ in the other direction than $\pi_2$. The
resulting space $Y$ is smooth.
A schematic picture of $Y$ is drawn in Figure~\ref{sfig:bu4}.

If $f:\widehat{X}\rightarrow Y$ denotes the contraction,
the formula~\eqref{eq:discrepancy} gives
$$
  K_{\widehat{X}}=f^*K_Y+E.
$$ 

We then use Proposition~\ref{prop:formula} to each component $E_j$,
taking $n=2$ in the formulas given there. We first claim that
$$
f^*f_*M=M+E,\qquad f^*f_*D=D+2E.
$$ 
We explain how to get these formulas from Figure~\ref{fig:comb-A4},
since this is the method we will use for more complicated
arrangements, but the reader may also want to glance at
Figure~\ref{fig:blowup}.

The first formula comes counting the planes in the projectivized
arrangement that intersect each given line $\pi(E_j)$ away from the
points blown-up. Indeed, these correspond to the components of $f_*M$
that contain $f(E_j)$.

Recall that $\pi(E_j)$ is an element in the $G$-orbit of $L_{12}$. The
last column (incident vertices) in the $L_{12}$-row of
Figure~\ref{fig:comb-A4} indicates that $L_{12}$ contains three
1-dimensional strata, two in the $G$-orbit of $L_{123}$ and one in the
$G$-orbit of $L_{134}$. The first two correspond to points that get
blown-up, the third one to a transverse intersection in the $G$-orbit
of $L_1\cap L_{34}$, which is the same as the $G$-orbit of $L_4\cap
L_{12}$. This implies that, when pulling-back $f_*M$ under $f$, we
will pick-up every component $E_j$ precisely one, hence the announced
formula.

The formula for $f^*f_*D$ follows from the fact that there are two
components of $f_*D$ that contain each component $f(E_k)$ of
$f(E)$. Indeed, these correspond to the 1-dimensional mirror
intersections $L$ with $\kappa_L>1$ that are contained in $\pi(E_k)$;
the number of such 1-dimensional intersections can be found once again
in the $L_{12}$ row of Figure~\ref{fig:comb-A4} (it is indicated by
the $2\times L_{123}$ in the column for indicident vertices).

Now we get (for every $\nu,\delta\in\Q$),
\begin{eqnarray*}
  &\left(K_Y+\nu f_*M+\delta f_*D\right)^3=\\
  &\left(K_{\widehat{X}}-E+\nu(M+E)+\delta(D+2E)\right)^3=\\
  &\left(\pi^*K_X+2D+E-E+\nu(M+E)+\delta(D+2E)\right)^3.
\end{eqnarray*}
Note that
$$
  \pi^*\pi_*M=M+6D+3E,
$$ 
since there are six mirrors through each point blown-up (see the third
column in the $L_{123}$ row of Figure~\ref{fig:comb-A4}), and
three mirrors containing each line blown-up (see the third column in
the $L_{12}$ row of Figure~\ref{fig:comb-A4}).

Hence the above formula can be written as
\begin{eqnarray*}
   &\left(\pi^*(K_X+\nu \pi_*M)+ (2+\delta-6\nu)D - 2(\nu-\delta)E\right)^3.
\end{eqnarray*}

We then have
\begin{eqnarray*}
  &(K_Y+\Delta)^3=\left(\lambda \pi^*H + \alpha D + \beta E\right)^3
\end{eqnarray*}
where $H$ denotes the class of a plane in $\P^3$, and
\begin{eqnarray*}
  &\lambda=-4+10\nu,\quad 
  \alpha=2+\delta-6\nu,\quad 
  \beta=-2(\delta-\nu).
\end{eqnarray*}

Finally, we get
\begin{eqnarray*}
  &(K_Y+\Delta)^3=\\
  &\lambda^3 + 5\alpha^3 + 10\beta^3\cdot E^3 + 3 \cdot 10 \cdot \beta^2 (\lambda\pi^*H\cdot E^2 + \alpha D\cdot E^2)=\\
  &\lambda^3 + 5\alpha^3 + 20\beta^3 - 30 \beta^2 ( \lambda + 2 \alpha ).
\end{eqnarray*}
To explain the last two equalities, the key point is that for each
irreducible component $E_j$ of $E$,
$$
  E_j|_{E_j}=-l_1-l_2,
$$
where $l_1$ and $l_2$ are the respective fibers in $\P^1\times\P^1$.

Note that when developing the cube, most cross-terms disappear because
$$
(\pi^*H)^2\cdot D = (\pi^*H)^2\cdot E = \pi^*H\cdot D^2 = \pi^*H\cdot D\cdot E = 0.
$$ 
Indeed, $H$ can be represented by a hyperplane not going through
any of the points blown-up, and $H^2$ can be represented by a line
that does not intersect any of the lines whose strict transform gets
blown-up.

Moreover, whenever $D_k$ intersects $E_j$, $D_k|_{E_j}=l_1$ (see the
discussion on p.~\pageref{prop:formula}), so $D_k^2\cdot E_k=0$. Also, we
can represent $H$ by a plane that is transverse to $\pi(E_j)$, so that
$\pi^*H|_{E_j}=l_1$, so we have
$$
D_k\cdot E_j^2=-1, \quad \pi^*H\cdot E_j^2=-1.
$$

Finally we have $D_j^3=1$ and
$$
E_j^3=E_j|_{E_j}\cdot E_j|_{E_j}=(-l_1-l_2)^2=2,
$$
see again the computations on p.~\pageref{prop:formula}.

The log-canonical divisor for $\CHL(A_4,8)$ is given as above for $\nu=1-\frac{1}{8}$, $\delta=2-\kappa_{L_{123}}=1-\frac{1}{2}$ (see the tables in
Figure~\ref{fig:comb-B4} on p.~\pageref{fig:comb-B4}). This gives
$$
(K_Y+\Delta)^3=\frac{33}{8},
$$
and finally
$$
\chi^{orb}(\CHL(A_4,8))=-\frac{11}{5120}.
$$
This agrees with the formula in~\cite{mcmullengaussbonnet}.



\subsection{The groups derived from the $G_{28}$ arrangement} \label{sec:g28}

The combinatorial properties of the $G_{28}$ arrangement are given in
the tables of Figure~\ref{fig:comb-g28}, p.~\pageref{fig:comb-g28}.

There are two orbits of mirrors of complex reflections in $G=G_{28}$,
each containing 12 hyperplanes. The $G$-invariant weight assignments
are parametrized by a pair $(p_1,p_2)$ of integers.

In fact there is an outer automorphism of $G$ exchanging the two
conjugacy classes of complex reflections, so the groups
$\CHL(p_1,p_2)$ and $\CHL(p_2,p_1)$ are isomorphic. Without loss of
generality, we may and will assume that $p_1\leqslant p_2$.

The Schwarz condition is of course satisfied for $(p_1,p_2)=(2,2)$, in
which case the group $\CHL(G_{28},2,2)$ is simply $G_{28}$. It is also
satisfied for $(p_1,p_2)=(2,3)$, which gives a parabolic group,
i.e. the signature of the invariant Hermitian form is $(3,0)$.

There are 11 other pairs $(p_1,p_2)$ with $p_1\leqslant p_2$ such that
where the Schwarz condition holds, listed in the table of
Figure~\ref{fig:comb-g28} (p.~\pageref{fig:comb-g28}). We will compute
volumes for all cases, grouping them in families where the blow-up
$\widehat{X}$ and the contracted space $Y$ have the same
description, hence the corresponding volume formulae are similar.

\subsubsection{The $G_{28}$ cases where no blow-up is needed} \label{sec:g28-noblowup}

There are two such groups, given by $(p_1,p_2)=(2,4)$ or $(3,3)$. As
above, we write $\nu_j=1-\frac{2}{p_j}$. For these two cases we have
$\widehat{X}=X=Y$, and $K_X+D$ is numerically equivalent to
$(-4+12\nu_1+12\nu_2)H$, so
$$
  (c_1^{orb}(\B^3/\CHL(G_{28},p_1,p_2))^3=\frac{1}{576}(-4+12\nu_1+\nu_2)^3,
$$
where the denominator 576=1152/2 is the order of the projective group
$\P(G_{28})$, in other words $|G_{28}|/|Z(G_{28})|$.

This gives
$$
  \chi^{orb}((\B^3/\CHL(G_{28},2,4)=-\frac{1}{1152},
$$
and
$$
  \chi^{orb}((\B^3/\CHL(G_{28},3,3)=-\frac{1}{144}.
$$

\subsubsection{The $G_{28}$ cases where we blow up points} \label{sec:g28-points}

There are 5 pairs of weights where we only blow up points. For
$(p_1,p_2)=(2,5)$ and $(2,6)$, we blow up the 12 points in the
$G$-orbit of $\P L_{234}$.  For $(p_1,p_2)=(3,4)$, $(3,6)$, $(4,4)$,
we blow up two $G$-orbits of points in $X=\P^3$, namely the 12 points
in the $G$-orbit of $\P L_{234}$ and the 12 points in the $G$-orbit of
$\P L_{123}$.

We treat the cases where we blow up two $G$-orbits of points in some
detail, the other ones (where we blow-up only one $G$-orbit) are
easier. Denote by $\pi:\widehat{X}\rightarrow X$ the corresponding
blow up, and by $M_1$ and $M_2$ the proper transform of the two orbits
of mirrors in $G_{28}$.

Write $D_1$ (resp. $D_2$) for the exceptional divisor above the $G$-orbit
of $\P L_{123}$ (resp. $\P L_{234}$). Note that the divisors $D_j$ both have
12 disjoint components.  We have
\begin{eqnarray*}
  & K_{\widehat{X}}=\pi^*K_X+2D_1+2D_2\\
  & \pi^*\pi_*M_1=6D_1+3D_2\\
  & \pi^*\pi_*M_2=3D_1+6D_2\\
\end{eqnarray*}
The last two formulae follow from the count of mirrors of each type
through $L_{123}$ (resp. $L_{234}$), see Figure~\ref{fig:comb-g28}
(for $j=1,2$, mirrors of type $j$ are those in the $G$-orbit of the
mirror $L_j$ of $r_j$). Note that $L_{123}$ is on 6 mirrors of type 1 and 3
mirrors of type 2, as is indicated by 6+3 in the third column of the
row headed $L_{123}$. Similarly, $L_{234}$ is on 3 mirrors of type 1
and 6 mirrors of type 2.

The relevant divisor for the orbifold pair is
$$
  \Delta = \nu_1 M_1 + \nu_2 M_2 + \delta_1 D_1 + \delta_2 D_2,
$$ 
where, $\nu_j=1-\frac{2}{p_j}$ and $\delta_1=2-\kappa_{L_{123}}$,
  $\delta_2=2-\kappa_{L_{234}}$.

We need to compute
\begin{eqnarray}
&(K_{\widehat{X}}+\Delta)^3\\
&=\left(\pi^*K_X + \nu_1 M_1 + \nu_2 M_2 + (2+\delta_1) D_1 + (2+\delta_2) D_2\right)\\
&=\left(\lambda\pi^*H + \alpha_1 D_1 + \alpha_2 D_2\right)^3\\
&=\lambda^3 + 12\alpha_1^3 + 12 \alpha_2^3 \label{eq:last}
\end{eqnarray}
where
\begin{eqnarray*}
  &\lambda = -4 + 12\nu_1 + 12\nu_2,\quad \alpha_1 = 2+\delta_1 - 6 \nu_1 - 3 \nu_2,\quad \alpha_2 = 2+\delta_2 - 3 \nu_1 - 6 \nu_2.
\end{eqnarray*}
The factors of 12 in equation~\eqref{eq:last} come from the fact that
each $D_j$, $j=1,2$ has 12 components, note also that $D_j^3=1$. When
developing the cube, the cross-terms do not contribute since the 24
components of $D_1+D_2$ are pairwise disjoint, and $H$ can be
represented by a plane not containing any of the 24 points blown-up.

Formula~\eqref{eq:last} gives
$$
\chi^{orb}(\B^3/\CHL(G_{28},3,4))=-\frac{23}{1152},
$$
$$
\chi^{orb}(\B^3/\CHL(G_{28},3,6))=-\frac{1}{36},
$$
$$
\chi^{orb}(\B^3/\CHL(G_{28},4,4))=-\frac{5}{144}.
$$

For the cases $(p_1,p_2)=(2,5)$ and $(2,6)$, the formula is the same,
except one removes the term corresponding to $D_1$ (i.e. the
exceptional above the $G$-orbit of $\P L_{123}$, which is not supposed
to get blown-up since $\kappa_{L_{123}}<1$).  In other words, with the
same notation for $\alpha$ and $\lambda$, we have
$$
  (K_{\widehat{X}}+\Delta)^3=\lambda^3 + 12\alpha_2^3.
$$

This gives
$$
\chi^{orb}(\B^3/\CHL(G_{28},2,5))=-\frac{13}{4500},
$$
$$
\chi^{orb}(\B^3/\CHL(G_{28},2,6))=-\frac{5}{1296}.
$$

In sections~\ref{sec:g28-lines1} through~\ref{sec:g28-lines3}, we
treat the $G_{28}$ cases where we need to blow up both points and
lines.

\subsubsection{The case $\CHL(G_{28},p_1,p_2)$ for $(p_1,p_2)=(2,8)$ or $(2,12)$} \label{sec:g28-lines1}

In this case we blow up the 12 points in the $G$-orbit of $\P
L_{234}$, and then the strict transform of the 16 lines in the
$G$-orbit of $\P L_{34}$. As before, we denote the corresponding
composition of blow-ups by $\pi:\widehat{X}\rightarrow X=\P^3$ and by
$f:\widehat{X}\rightarrow Y$ the relevant contraction, see
section~\ref{sec:chl}.

Since each copy of $L_{34}$ contains 3 copies of $L_{234}$,
$f:\widehat{X}\rightarrow Y$ is crepant, i.e. $f^*K_Y=K_X$ (see
equation~\eqref{eq:discrepancy} for $n=3$).

On $\widehat{X}$ we have $D$ (exceptional divisor with 12 components,
above the $G$-orbit of $\P L_{234}$), $E$ (exceptional divisor with 16
components, above the $G$-orbit of $\P L_{34}$), $M_1$ and $M_2$
(strict transform of the $G$-orbit of mirrors of reflections in $G$,
both have 12 components). We will need the following formulae:
\begin{eqnarray*}
  &K_{\widehat{X}}=f^*K_Y\\
  &f^*f_*M_1=M_1+\frac{3}{2}E,\qquad f^*f_*M_2=M_2\\
  &f^*f_*D=D+\frac{3}{2}E\\
  &K_{\widehat{X}}=\pi^*K_X+2D+E\\
  &\pi^*\pi_*M_1=M_1+3D,\qquad \pi^*\pi_*M_2=M_2+6D+3E
\end{eqnarray*}

In the formula for $f^*f_*M_1$, the denominator 2 comes from part~(2)
of Proposition~\ref{prop:formula}, and the numerator 3 comes from the
count of the number of components of $f_*M_1$ that contain the image
of a component of $E$ in $Y$. The latter number is given by the number
of mirrors in the $G$-orbit of $L_1$ that intersect $L_{34}$
transversely, away from the $G$-orbit of $L_{234}$ (these are the
points that get blown-up to get $\widehat{X}$). It is indicated by the
occurrence of $3\times L_{134}$ in the last column of the row of
Figure~\ref{fig:comb-g28} headed $L_{34}$.

Similarly, in the formula for $f^*f_*D$, the denominator 2 comes from
formula in Proposition~\ref{prop:formula}(3), and the numerator 3
comes from the number of points of $\P L_{34}$ that are in the
$G$-orbit of $\P L_{234}$ (see the $3\times L_{234}$ in the row of
Figure~\ref{fig:comb-g28} headed $L_{34}$).

The formulae for $\pi^*\pi_*M_j$ ($j=1,2$) follow from the count of
mirrors of each type containing $L_{234}$ (the 3+6 in the table
indicates that it is contained in 3 mirrors of type 1, and 6 mirrors
of type 2) and $L_{34}$ (0+3 indicates that it is contained in 3
mirrors of type 2).

The same computations as in section~\ref{sec:a4-lines} then give
\begin{equation}\label{eq:g28-cube}
  \left(K_X+\Delta\right)^3=\left(\lambda\pi^*H+\alpha D+\beta E\right)^3
\end{equation}
where
\begin{eqnarray*}
  &\lambda = -4 + 12\nu_1 + 12\nu_2\\
  &\alpha = 2+\delta - 3\nu_1 - 6\nu_2\\
  &\beta = 1 + \frac{3}{2}\nu_1 - 3\nu_2 + \frac{3}{2}\delta\\
\end{eqnarray*}
and we need to take $\nu_j=1-\frac{2}{p_j}$, and $\delta=2-\kappa_{L_{234}}$.

When developing the cube in equation~\eqref{eq:g28-cube}, most
cross-term disappear for the same reason as in
section~\ref{sec:a4-lines}. Once again only terms of the form
$D_k\cdot {E_j}^2$ or $\pi^*H\cdot {E_j}^2$ remain, where the $D_k$
(resp. $E_j$) denotes the $k$-th component of $D$ (resp. the $j$-th
component of $E$).

Recall that $E_j|_{E_j}=-l_1-2l_2$, where $l_1$ is the class in
$Pic(E_j)$ that projects to a line in $\P^3$ (see
equation~\eqref{eq:self}). Moreover, $\pi^*H$ restricts to $l_2$, and
$D_k$ restricts to either 0 or $l_2$ (depending on whether $D_k$ and
$E_j$ intersect at all), see the discussion
on p.~\pageref{prop:formula}. This gives $\pi^*H\cdot E_j^2=-1$, and
$D_k\cdot E_j^2=-1$ (or 0 if $D_k$ and $E_j$ are disjoint).

Note that $E$ has 16 components, and $\pi^*H\cdot E_j^2=-1$ for each
$j$, we have $\pi^*H\cdot E^2=-16$. Also each component $E_j$ of $E$
intersects precisely 3 components of $D$ (see the occurrence of
$3\times L_{234}$ in the row for $L_{34}$ of
Figure~\ref{fig:comb-g28}), and $D_k\cdot E_j^2=-1$ for each $j$, so
$D\cdot E^2=-16\cdot 3$.

Finally, we get
$$
\lambda^3 + 12\alpha^3 + 16\cdot 4\cdot \beta^3 - 3\cdot 16\cdot \lambda\beta^2 - 3\cdot 16\cdot 3\cdot \alpha\beta^2.
$$
For $(p_1,p_2)=(2,8)$, $\delta=2-\kappa_{L_{234}}=1-\frac{1}{2}$, this gives 
$$
  \chi^{orb}(\CHL(G_{28},2,8))=-\frac{11}{3072}.
$$ 
For $(p_1,p_2)=(2,12)$, $\delta=2-\kappa_{L_{234}}=1-\frac{2}{3}$, this gives 
$$
  \chi^{orb}(\CHL(G_{28},2,12))=-\frac{23}{10368}.
$$ 

\subsubsection{The case $\CHL(G_{28},6,6)$} \label{sec:g28-lines2}

In this case we blow up 12 points in the $G$-orbit of $\P L_{123}$,
the 12 points in the $G$-orbit of $\P L_{234}$, and then
the strict transform of the 18 lines in the $G$-orbit of $\P L_{23}$.

Since each copy of $L_{23}$ contains 2 copies of $L_{234}$ and 2
copies of $L_{123}$ (see the occurrence of $2\times L_{123}$ and
$2\times L_{234}$ in the column for ``Incident vertices'' of the row
headed $L_{23}$ in Figure~\ref{fig:comb-g28}), we have
$$
  K_{\widehat{X}}=f^*K_Y-\frac{n-3}{n-1}E
$$
with $n=2+2=4$, i.e. $f^*K_Y=K_{\widehat{X}}+\frac{1}{3}E$.

On $\widehat{X}$ we have $D_1,D_2,E,M_1,M_2$; here $D_1$
(resp. $D_2$) is the exceptional divisor above the $G$-orbit of $\P
L_{123}$ (resp. $\P L_{234}$). We have the following:
\begin{eqnarray*}
  &K_{\widehat{X}}=f^*K_Y+\frac{1}{3}E_1\\
  &f^*f_*M_1=M_1,\qquad f^*f_*M_2=M_2\\
  &f^*f_*D_1=D_1+\frac{2}{3}E,\qquad f^*f_*D_2=D_2+\frac{2}{3}E\\
  &K_{\widehat{X}}=\pi^*K_X+2D_1+2D_2+E\\
  &\pi^*\pi_*M_1=M_1+6D_1+3D_2+2E,\qquad \pi^*\pi_*M_2=M_2+3D_1+6D_2+2E
\end{eqnarray*}
The claim about $f^*f_*M_j$ follows from the fact that no mirror
intersects $L_{23}$ transversely away from the $G$-orbit of $L_{123}$
and away from the $G$-orbit of $L_{234}$.

The claim about $f^*f_*D_j$ follows from
Proposition~\ref{prop:formula}(3), and the fact that each $E_k$
contains two points of the $G$-orbit of $\P L_{123}$, and two points
in the $G$-orbit of $\P L_{234}$.

Computations similar to those in section~\ref{sec:g28-lines1} now give
$$
  \left(K_X+\Delta\right)^3=\left(\lambda\pi^*H+\alpha D+\beta E\right)^3
$$
where
\begin{eqnarray*}
  &\alpha_1 = 2 + \delta_1 - 6\nu_1 - 3\nu_2,\quad
  \alpha_2 = 2 + \delta_2 - 3\nu_1 - 6\nu_2\\
  &\lambda = -4 + 12\nu_1 + 12\nu_2, \quad \beta = \frac{4}{3} - 2\nu_1 - 2\nu_2 + \frac{2}{3}\delta_1 + \frac{2}{3}\delta_2.\\
\end{eqnarray*}
In the above formulae, we now take $\nu_j=1-\frac{2}{p_j}$,
$\delta_1=2-\kappa_{L_{123}}$, $\delta_2=2-\kappa_{L_{234}}$.

Recall that we have $E_j^3=2\cdot 4-2=6$ for every component $E_j$ of
$E$, see equation~\eqref{eq:self}. The same analysis of the
cross-terms as in section~\ref{sec:g28-lines1} gives
$$
\left(K_X+\Delta\right)^3=
\lambda^3 + 12\alpha_1^3 + 12\alpha_2^3 + 18\cdot 6\cdot \beta^3 
- 3 \cdot 18\cdot \lambda\beta^2 
- 3 \cdot 18\cdot 2\cdot \alpha_1\beta^2 
- 3 \cdot 18\cdot 2\cdot \alpha_2\beta^2 
$$
which for $(p_1,p_2)=(6,6)$, $\delta_1=\delta_2=0$ gives 
$$
  \chi^{orb}(\CHL(G_{28},6,6))=-\frac{5}{144}.
$$ 

\subsubsection{The case $\CHL(G_{28},3,12)$} \label{sec:g28-lines3}

In this case we blow up the 12 points in the $G$-orbit of $\P
L_{123}$, the 12 points in the $G$-orbit of $\P L_{234}$, then the
strict transform of the 18 lines in the $G$-orbit of $\P L_{23}$, and
the strict transform of the 16 lines in the $G$-orbit of
$\P L_{34}$.

On $\widehat{X}$ we now have $D_1,D_2,E_1,E_2,M_1,M_2$ (where $D_1$
corresponds to $L_{123}$, $D_2$ to $L_{234}$, $E_1$ to $L_{23}$, $E_2$
to $L_{34}$).
\begin{eqnarray*}
  &K_{\widehat{X}}=f^*K_Y-\frac{1}{3}E_1\\
  &f^*f_*M_1=M_1+\frac{3}{2}E_2,\qquad f^*f_*M_2=M_2\\
  &f^*f_*D_1=D_1+\frac{2}{3}E_1,\qquad f^*f_*D_2=D_2+\frac{2}{3}E_1+\frac{3}{2}E_2\\
  &K_{\widehat{X}}=\pi^*K_X+2D_1+2D_2+E_1+E_2\\
  &\pi^*\pi_*M_1=M_1+6D_1+3D_2+2E_1,\qquad \pi^*\pi_*M_2=M_2+3D_1+6D_2+2E_1+3E_2
\end{eqnarray*}
The same computations as before now give
$$
  \left(K_X+\Delta\right)^3=\left(\lambda\pi^*H+\alpha D+\beta E\right)^3
$$
where
\begin{eqnarray*}
  &\lambda = -4 + 12\nu_1 + 12\nu_2\\
  &\alpha_1 = 2+\delta_1 - 6\nu_1 - 3\nu_2,\qquad
  \alpha_2 = 2+\delta_2 - 3\nu_1 - 6\nu_2\\
  &\beta_1 = \frac{4}{3} - 2\nu_1 - 2\nu_2 + \frac{2}{3}\delta_1 + \frac{2}{3}\delta_2, \quad \beta_2 = 1 + \frac{3}{2}\nu_1 - 3\nu_2 + \frac{3}{2}\delta_1\\
\end{eqnarray*}
Developing the cube, we get
{\scriptsize
\begin{eqnarray*}
&&\lambda^3 + 12\alpha_1^3 + 12\alpha_2^3
          + 18\beta_1^3\cdot 9 + 16\beta_2^3\cdot 4
          + 3 \left(\pi^*H\cdot E_1^2
          + \pi^*H\cdot E_2^2
          + D_1\cdot E_1^2
          + D_2\cdot E_1^2
          + D_2\cdot E_2^2\right)\\
&&= \lambda^3 + 12\alpha_1^3 + 12\alpha_2^3 + 162\beta_1^3 + 64\beta_2^3
- 3 \left(  
\cdot 18\cdot \lambda\beta_1^2 
+ \cdot 16\cdot \lambda\beta_2^2 
+  \cdot 18\cdot 2\cdot \alpha_1\beta_1^2 
+ \cdot 18\cdot 2\cdot \alpha_2\beta_1^2 
+  \cdot 16\cdot 3\cdot \alpha_2\beta_2^2\right).
\end{eqnarray*}
}
For $(p_1,p_2)=(3,12)$, $\delta_1=1-\frac{1}{2}$, $\delta_2=1-1=0$, we get
$$
  \chi^{orb}(\CHL(G_{28},3,12))=-\frac{23}{1152}.
$$

\subsection{The groups derived from the $G_{29}$ arrangement} \label{sec:g29}

The combinatorial properties of the $G_{29}$ arrangement are given in
the tables of Figure~\ref{fig:comb-g29}, p.~\pageref{fig:comb-g29}.

Here the group has a single orbit of mirrors of reflections, so the
corresponding lattices $\CHL(G_{29},p)$ are indexed by a single
integer $p$. The hyperbolic cases that satisfy the Schwarz condition
correspond to $p=3$ or $4$.

The volume computations are similar to the ones in
section~\ref{sec:a4-points} or section~\ref{sec:g28-points}, since we
only blow up points, i.e. in the notation of
section~\ref{sec:volumes}, $\widehat{X}=Y$.

For $p=3$, we need to blow up the $G$-orbit of $L_{234}$, since
$\kappa_{L_{234}}=1+\frac{1}{3}>1$. In this case, we get
$$
  \left(K_{\widehat{X}}+\Delta\right)^3=\lambda^3+20\alpha^3,
$$
where
\begin{eqnarray*}
  & \lambda=-4+40\nu\\
  & \alpha=2+\delta-12\nu,
\end{eqnarray*}
where $\nu=1-\frac{2}{p}$ and $\delta=2-\kappa_{234}=1-\frac{1}{3}$.
This gives
$$
\chi^{orb}(\CHL(G_{29},3))=-\frac{323}{12960}.
$$

For $p=4$, we need to blow up the $G$-orbit of $L_{234}$ (which gives
20 points in $X=\P^3$) and the $G$-orbit of $L_{124}$ (which gives 40
points in $X=\P^3$), see Figure~\ref{fig:comb-g29}. The formula is
similar to the one for $p=3$, we get
$$
\left(K_{\widehat{X}}+\Delta\right)^3=\lambda^3+20\alpha_1^3+40\alpha_2^3,
$$
where
\begin{eqnarray*}
  & \lambda = -4 + 40\nu\\
  & \alpha_1 = 2 + \delta_1 - 12\nu\\
  & \alpha_2 = 2 + \delta_2 - 9\nu.
\end{eqnarray*}
Taking $\nu=1-\frac{1}{4}$, $\delta_1=2-\kappa_{234}=0$,
$\delta_2=2-\kappa_{124}=1-\frac{1}{2}$, we get
$$
\chi^{orb}(\CHL(G_{29},4))=-\frac{13}{160}.
$$

\subsection{The groups derived from the $G_{30}$ arrangement} \label{sec:g30}

The combinatorial properties of the $G_{30}$ arrangement are given in the tables of
Figure~\ref{fig:comb-g30}, p.~\pageref{fig:comb-g30}. 

Again, the group has a single orbit of mirrors of reflections. The
Schwarz condition holds (and the group is hyperbolic) for
$\CHL(G_{30},p)$, $p=3$ or $5$.

For $p=3$, we only blow up points (given by the 60 points
corresponding to the $G$-orbit of $L_{234}$), so the computation is
similar to the one in section~\ref{sec:a4-points}. We get
$$
\left(K_{\widehat{X}}+\Delta\right)^3=\lambda^3+60\alpha^3,
$$
where
\begin{eqnarray*}
  & \lambda=-4+60\nu\\
  & \alpha=2+\delta-15\nu.
\end{eqnarray*}
Taking $\nu=1-\frac{2}{3}$, $\delta=1-\frac{2}{3}$, we get
$$
\chi^{orb}(\CHL(G_{30},3))=-\frac{52}{2025}.
$$

For $p=5$, we blow up 300 points corresponding to the $G$-orbit of
$L_{123}$, the 60 points corresponding to the $G$-orbit of $L_{234}$,
and then the strict transform of the 72 lines corresponding to the
$G$-orbit of $L_{34}$.

We denote the corresponding exceptional divisors by $D_1,D_2,E$, and
note
\begin{eqnarray*}
  &K_{\widehat{X}}=f^*K_Y-\frac{1}{2}E\\
  &f^*f_*L=L+\frac{5}{4}E\\
  &f^*f_*D_1=D_1,\qquad f^*f_*D_2=D_2+\frac{5}{4}E\\
  &K_{\widehat{X}}=\pi^*K_X+2D_1+2D_2+E\\
  &\pi^*\pi_*L=L+6D_1+15D_2+5E.
\end{eqnarray*}
The same computations as before now give
$$
  \left(K_X+\Delta\right)^3=\left(\lambda\pi^*H+\alpha_1 D_1+\alpha_2 D_2+\beta E\right)^3
$$
where
\begin{eqnarray*}
  &\lambda=-4+60\nu\\
  &\alpha_1=2+\delta-1-6\nu,\qquad
  \alpha_2=2+\delta_2-15\nu\\
  &\beta=\frac{3}{2}-\frac{15}{4}\nu+\frac{5}{4}\delta_2.
\end{eqnarray*}
Using the combinatorics and the self intersection of $E_1$ and $E_2$
(see the previous sections), we get
$$
\lambda^3 + 300\alpha_1^3 + 60\alpha_2^3 + 72\cdot 8\cdot \beta^3 
- 3 \cdot 72\cdot \lambda\beta^2 
- 3 \cdot 72\cdot 5\cdot \alpha_2\beta^2. 
$$
We then take $\nu=1-\frac{2}{5}$, $\delta_1=1-\frac{1}{5}$, $\delta_2=1-2=-1$, and get
$$
  \chi^{orb}(\CHL(G_{30},5))=-\frac{41}{1125}.
$$

\subsection{The groups derived from the $G_{31}$ arrangement} \label{sec:g31}

The combinatorial properties of the $G_{31}$ arrangement are given in
the tables of Figure~\ref{fig:comb-g31}, p~\pageref{fig:comb-g31}.

There are two values of $p$ such that the Schwarz condition and the
group $\CHL(G_{31},p)$ is hyperbolic, namely $p=3$ or $p=5$.

For $p=3$, we need to blow up the 60 points corresponding to the
$G$-orbit of $L_{125}$. We get
$$
  \left(K_{\widehat{X}}+\Delta\right)^3=\lambda^3+60\alpha^3,
$$
where
\begin{eqnarray*}
  & \lambda=-4+60\nu\\
  & \alpha=2+\delta-15\nu.
\end{eqnarray*}
Taking $\nu=1-\frac{2}{3}$, $\delta=1-\frac{2}{3}$, we get
$$
\chi^{orb}(\CHL(G_{31},3))=-\frac{13}{810}.
$$

For $p=5$, we blow up the 60 points corresponding to the $G$-orbit of
$L_{125}$, the 480 points corresponding to the $G$-orbit of $L_{235}$,
and then the 30 lines corresponding to the $G$-orbit of $L_{14}$.

We denote the corresponding exceptionals by $D_1,D_2,E$, and note
\begin{eqnarray*}
  &K_{\widehat{X}}=f^*K_Y-\frac{3}{5}E\\
  &f^*f_*L=L\\
  &f^*f_*D_1=D_1+\frac{6}{5}E,\qquad f^*f_*D_2=D_2\\
  &K_{\widehat{X}}=\pi^*K_X+2D_1+2D_2+E\\
  &\pi^*\pi_*L=L+15D_1+6D_2+6E.
\end{eqnarray*}
The same computations as before now give
$$
  \left(K_X+\Delta\right)^3=\left(\lambda\pi^*H+\alpha_1 D_1+\alpha_2 D_2+\beta E\right)^3
$$
where
\begin{eqnarray*}
  &\lambda = -4+60\nu\\
  &\alpha_1 = 2+\delta_1 - 15\nu\\
  &\alpha_2 = 2+\delta_2 - 6\nu\\
  &\beta = \frac{8}{5} - 6\nu + \frac{6}{5}\delta_2.
\end{eqnarray*}
Using the combinatorics and the self intersection of $E_1$ and $E_2$
(see the previous sections), we get
$$
\lambda^3 + 60\alpha_1^3 + 480\alpha_2^3 + 30\cdot 10\cdot \beta^3 
- 3 \cdot 30\cdot \lambda\beta^2 
- 3 \cdot 30\cdot 6\cdot \alpha_1\beta^2. 
$$
Finally, taking $\nu=1-\frac{2}{5}$, $\delta_1=-1$, $\delta_2=1-\frac{1}{5}$, we get
$$
  \chi^{orb}(\CHL(G_{31},5))=-\frac{41}{1125}.
$$

\subsection{The groups derived from the $B_4$ arrangement} \label{sec:b4}


For completeness, we compute the volumes of the CHL groups associated
to the $B_4$ arrangement, even though the corresponding volumes are
known. Indeed, the lattices of the form $\CHL(B_4,p_1,p_2)$ are
commensurable to certain Deligne-Mostow groups (see
section~\ref{sec:dm}).

The combinatorial properties of the $B_4$ arrangement are given in the
tables of Figure~\ref{fig:comb-B4}, p.~\pageref{fig:comb-B4}. In this
case, the group $G=W(B_4)$ has two orbits of mirrors of complex
reflections. In the numbering used in Figure~\ref{fig:comb-B4}, the
mirror of $r_1$ is not in the same orbit as the mirror of $r_2$ (but
the mirror of $r_3$ is in the same orbit as the mirror of $r_2$, since
the braid relation $r_2r_3r_2=r_3r_2r_3$ implies that $r_2$ and $r_3$
are conjugate in $W(B_4)$, since $r_3=r_2r_3\cdot r_2\cdot
(r_2r_3)^{-1}$).

The $G$-invariant weight assignments are determined by the two weights
$\kappa_1=1-\frac{2}{p_1}$ and $\kappa_2=1-\frac{2}{p_2}$, where $p_j$
are integers $\geq 2$. As before, we denote by $\CHL(B_4,p_1,p_2)$ the
corresponding group.

For $(p_1,p_2)=(n,2)$ for some $n\geq 2$, the group
$\CHL(B_4,p_1,p_2)$ turns out to be finite (in the Shephard-Todd
notation, it is given by the group $G(n,1,4)$, which is imprimitive
for $n>2$). For $(p_1,p_2)=(2,3)$, the Hermitian form preserved by the
group is degenerate of signature (3,0), and the corresponding group
gives a complex affine crystallographic group acting on $\C^3$ (see
section~5 of~\cite{chl}).

The other pairs $(p_1,p_2)$ where the Schwarz condition holds are all
hyperbolic (i.e. the group $\CHL(B_4,p_1,p_2)$ preserves a Hermitian
form of signature $(3,1)$). The list of these pairs is given in
Figure~\ref{fig:comb-B4}. We treat them separately over
sections~\ref{sec:B4_1}--\ref{sec:B4_7}, according to the dimension of
the strata of the arrangement that need to be blown-up in order to
describe the quotient $\B^3/\CHL(B_4,p_1,p_2)$.


\subsubsection{The groups $\CHL(B_4,2,4)$ and $\CHL(B_4,3,3)$} \label{sec:B4_1}

In these cases, no blow-up is needed, since $\kappa_L\leq 1$ for every
irreducible mirror intersection in the arrangement. In other words, in
the notation of section~\ref{sec:volumes}, we have
$X=\widehat{X}=Y$. Since there are 4 mirrors in the first orbit and 12
mirrors in the second orbit, the log-canonical divisor is given on the
level of $\P^3$ by 
$(-4+4\nu_1+12\nu_2)H$, where $\nu_j=1-\frac{2}{p_j}$ and $H$ denotes
the hyperplane class.

Up to removal of the cusp, the ball quotient is given by $X/G$ where
$G=W(B_4)$ (with a different orbifold structure than the one coming
from this finite quotient), so we have
$$
c_1^{orb}(\B^3/\CHL(B_4,p_1,p_2))=\frac{1}{|\P(G)|}(K_X+\Delta)^3=\frac{1}{192}\left(-4+4\nu_1+12\nu_2\right)^3,
$$
which gives
$$
\chi^{orb}(\B^3/\CHL(B_4,p_1,p_2))=-\frac{1}{16\cdot 192}\left(-4+4\nu_1+12\nu_2\right)^3.
$$ 

For $(p_1,p_2)=(3,3)$, we get $-1/1296$ , which is the orbifold Euler
characteristic of $\Gamma_{\mu,\Sigma}$ for hypergeometric exponents
$\mu=(1,1,1,1,3,5)/6$, $\Sigma\simeq S_4$
(see~\cite{mcmullengaussbonnet}).

For $(p_1,p_2)=(2,4)$, we get $-1/384$, which is the Euler
characteristic of $\Gamma_{\mu,\Sigma_0}$ for $\mu=(1,1,1,1,1,3)/4$ and
$\Sigma_0\subset \Sigma\simeq S_5$ fixing one of the 5 equal
weights. This is coherent with the value given
in~\cite{mcmullengaussbonnet}, which is $-1/1920$, since
$\Gamma_{\mu,\Sigma_0}$ has index 5 in $\Gamma_{\mu,\Sigma}$.


\subsubsection{The groups $\CHL(B_4,p_1,p_2)$ for $(p_1,p_2)=(3,4)$, $(4,3)$, $(4,4)$ and $(6,3)$} \label{sec:B4_2}

In these cases, there is a single $G$-orbit of irreducible mirror
intersections $L$ with $\kappa_L>1$, namely the $G$-orbit of $L_{123}$
(see Table~\ref{fig:comb-B4}). We then have $Y=\widehat{X}$, where $X$
is obtained from $X=\P^3$ by blowing up the $G$-orbit of $L=L_{123}$,
which gives 4 points in $\P^3$. The relevant log-canonical divisor has
the form
$$
  K_{\widehat{X}} 
  + \nu_1 M_1 
  + \nu_2 M_2 
  + \delta D,
$$ 
where $\nu_j=1-\frac{2}{p_j}$ and $\delta=2-\kappa_L$.
 
Note that $M_1$ has 4 components, whereas $M_2$ has 12, see
Figure~\ref{fig:comb-B4} on p.~\pageref{fig:comb-B4}.
Note also that $K_{\widehat{X}}=\pi^*K_X+2D$, and 
$$
  \pi^*\pi_*M_1 = M_1+3D, \quad \pi^*\pi_*M_2 = M_2+6D, 
$$
since $L_{123}$ is on 3 mirrors in the first orbit, and 6 mirrors in the second orbit.

Now the log-canonical divisor can be rewritten as
$$
  \pi^*(K_X+A)+\alpha D,
$$
where $A = K_X + \nu_1 M_1 + \nu_2 M_2$ and 
$$
  \alpha= 2 + \delta - 3\nu_1 - 6\nu_2.
$$

Finally, observe that $K_X+A$ is linearly equivalent to $\lambda H$ where 
$$
  \lambda=-4+4\nu_1+12\nu_2,
$$
and we get 
$$
   (c_1(\widehat{X},\mathcal{D}))^3=\lambda^3+4\alpha^3.
$$ 
Indeed, $D$ has 4 components (corresponding to the fact that the
$G$-orbit of $L_{123}$ has 4 points), and for each component $D_j$, we
have $D_j^3=1$.

For $(p_1,p_2)=(3,4)$, we get 
$$
   \chi^{orb}(\B^3/\CHL(B_4,3,4))=-\frac{31}{3456}.
$$ 
This is the same as the value of the orbifold Euler characteristic of
the Deligne-Mostow group $\Gamma_{\mu,\Sigma}$ for
$\mu=(3,3,3,3,5,7)/12$, $\Sigma=S_4$ (note that this is actually the
non-arithmetic lattice in $PU(3,1)$ constructed by Deligne and
Mostow).

For $(p_1,p_2)=(3,4)$, we get 
$$
   \chi^{orb}(\B^3/\CHL(B_4,4,3))=-\frac{23}{10368}.
$$
This is the same as the value of the orbifold Euler characteristic of
the Deligne-Mostow group $\Gamma_{\mu,\Sigma}$ for
$\mu=(2,2,2,2,7,9)/12$, $\Sigma=S_4$.

For $(p_1,p_2)=(4,4)$, we get 
$$
   \chi^{orb}(\B^3/\CHL(B_4,4,4))=-\frac{1}{96}.
$$
This is the same as the value of the orbifold Euler characteristic of
the Deligne-Mostow group $\Gamma_{\mu,\Sigma_0}$ for
$\mu=(1,1,1,1,2,2)/4$, $\Sigma_0=S_4$ (the maximal Deligne-Mostow
lattice for these weights corresponds to $\Sigma=S_4\times S_2$, and
it has orbifold Euler characteristic $-1/192$).

For $(p_1,p_2)=(6,3)$, we get 
$$
   \chi^{orb}(\B^3/\CHL(B_4,6,3))=-\frac{1}{324}.
$$
This is the same as the value of the orbifold Euler characteristic of
the Deligne-Mostow group $\Gamma_{\mu,\Sigma_0}$ for
$\mu=(1,1,1,1,4,4)/6$, $\Sigma_0=S_4$ (again, the maximal
Deligne-Mostow lattice for these weights corresponds to
$\Sigma=S_4\times S_2$, so its Euler characteristic is $-1/648$).


\subsubsection{The groups $\CHL(B_4,p_1,p_2)$ for $(p_1,p_2)=(2,5)$, $(2,6)$ and $(3,6)$} \label{sec:B4_3}

This case is similar to the previous one, except that now we need to
blow up the 4 points in the $G$-orbit of $L_{123}$, as well as the 8
points in the $G$-orbit of $L_{234}$.

The relevant log-canonical divisor has the form
$$
  K_{\widehat{X}} 
  + \nu_1 M_1 
  + \nu_2 M_2 
  + \delta_1 D_1
  + \delta_2 D_2.
$$ 
where $\nu_j=1-\frac{2}{p_j}$, and $\delta_1=2-\kappa_{L{123}}$,
$\delta_2=2-\kappa_{L{234}}$. Here we denote by $D_1$ (resp. $D_2$)
the exceptional divisors above the $G$-orbit of the projectivization
of $L_{123}$ (resp. $L_{234}$).

Note that $M_1$ has 4 components, whereas $M_2$ has 12, see
Figure~\ref{fig:comb-B4}.
Note also that $K_{\widehat{X}}=\pi^*K_X+2D_1+2D_2$, and 
$$
  \pi^*\pi_*M_1 = M_1+3D_1, \quad \pi^*\pi_*M_2 = M_2+6D_1+6D_2.
$$ 
Indeed, each element in the orbit of $L_{123}$ lies on 3 components of
$M_1$ and 6 components of $M_2$, and each element in the orbit of
$L_{234}$ lies on (no component of $M_1$ and) 6 components of $M_2$.

We get
$$
(K_{\widehat{X}}+\Delta)^3=(-4+4\nu_1+12\nu_2)^3 + 4 (2+\delta-3\nu_1)^3 + 8 (2+\delta_2-6\nu_1-6\nu_2)^3.
$$

For $(p_1,p_2)=(2,5)$, we take $\nu_j=1-2/p_j$,
$\delta_1=\delta_2=1-1/5$ (see the tables in
Figure~\ref{fig:comb-B4}), we get
$$
(c_1^{orb})^3=-\frac{1}{192}(K_{\widehat{X}}+\Delta)^3=-\frac{52}{375},
$$
so
$$
\chi^{orb}(\B^3/\CHL(B_4,2,5))=-\frac{13}{1500}.
$$
This agrees with the Euler characteristic of the Deligne-Mostow
lattice $\Gamma_{\mu,\Sigma}$ with $\mu=(2,3,3,3,3,6)/10$,
$\Sigma=S_4$.

For $(p_1,p_2)=(2,6)$, we take $\delta_1=\delta_2=1-1/3$ and get
$$
\chi^{orb}((\B^3/\CHL(B_4,2,6)))=-\frac{5}{432}.
$$
This agrees with the Euler characteristic of the Deligne-Mostow
lattice $\Gamma_{\mu,\Sigma}$ with $\mu=(1,2,2,2,2,3)/6$,
$\Sigma=S_4$.

For $(p_1,p_2)=(3,6)$, we take $\delta_1=\delta_2=1-1/3$ and get
$$
\chi^{orb}((\B^3/\CHL(B_4,3,6)))=-\frac{5}{432}.
$$
This agrees with the Euler characteristic of the Deligne-Mostow
lattice $\Gamma_{\mu,\Sigma_0}$ with $\mu=(1,1,1,1,1,1)/6$,
$\Sigma_0=S_4$. Note that this has index $30=6!/4!$ in the lattice
$\Gamma_{\mu,\Sigma}$ with $\Sigma=S_6$.


\subsubsection{The cases $\CHL(B_4,(p_1,p_2))$ with $(p_1,p_2)=(6,4)$ or $(12,3)$} \label{sec:b4-lines1} \label{sec:B4_4}

Here there are both lines and planes among the mirror intersections
$L$ that satisfy $\kappa_L>1$, which correspond to the $G$-orbit of
$L_{123}$ (this gives 4 points in $\P^3$) and the $G$-orbit of
$L_{12}$ (this gives 6 lines in $\P^3$).

We denote by $M_1$ and $M_2$ the strict tranform in $\widehat{X}$ of
the two $G$-orbits of mirrors ($M_1$ has 4 components, whereas $M_2$
has 12). We denote by $D$ (resp. $E$) the exceptional divisor in
$\widehat{X}$ above the $G$-orbit of $L_{123}$ (resp. the $G$-orbit of
$L_{12}$).

We need to compute
$$
\left(K_Y+\nu_1f_*M_1+\nu_2f_*M_2+\delta f_*D\right)^3,
$$
where $\nu_j=1-\frac{2}{p_2}$ and $\delta=2-\kappa_{L_{123}}$.

We have $K_{\widehat{X}}=f^*K_Y+E$ (see equation~\eqref{eq:discrepancy}), and
one checks using the combinatorics of the arrangement that
$$
f^*f_*M_1=M_1,\quad f^*f_*M_2=M_2+2E,\quad f^*f_*D=D+2E.
$$
Note also that
$$
\pi^*\pi_*M_1=M_1+3D+2E,\quad 
\pi^*\pi_*M_2=M_2+6D+2E,
$$ 
because for each $j$, $\pi_*D_j$ is on 3 mirrors in the first orbit,
and 6 mirrors in the second orbit, and $\pi^*E_j$ lies on 2 mirrors
from each orbit.

This gives
\begin{eqnarray*}
&(K_Y+\Delta)^3=\left(\lambda\pi^*H+\alpha D+\beta E\right)^3
\end{eqnarray*}
where
\begin{eqnarray*}
  &\lambda=-4+4\nu_1+12\nu_2,\quad
  \alpha=2+\delta-3\nu_1-6\nu_2,\quad
  \beta=2(\delta-\nu_1).
\end{eqnarray*}
Finally we get
\begin{eqnarray*}
  &(K_Y+\Delta)^3= \lambda^3 + 4\alpha^3 + 6\beta^3\cdot 2 
  + 3\lambda\beta^2 \pi^*H\cdot E^2 
  + 3\alpha\beta^2  D\cdot E^2\\
  &= \lambda^3 + 4\alpha^3 + 6\beta^3\cdot 2 
  - 3 \cdot 6\cdot \lambda\beta^2 
  - 3\cdot 6\cdot 2\cdot \alpha\beta^2.
\end{eqnarray*}

For $p_1=6$, $p_2=4$, we take $\nu_j=1-\frac{2}{p_j}$ and $\delta=1-\frac{2}{3}$, this gives
$$
\chi^{orb}(\CHL(B_4,6,4))=-\frac{31}{3456},
$$
and for $p_1=12$, $p_2=3$, we take $\nu_j=1-\frac{2}{p_j}$ and $\delta=1-\frac{1}{2}$, this gives
$$
\chi^{orb}(\CHL(B_4,12,3))=-\frac{23}{10368},
$$ 
as it should in comparison with the values expected
from~\cite{mcmullengaussbonnet}.


\subsubsection{The cases $\CHL(B_4,(p_1,p_2))$ with $(p_1,p_2)=(6,6)$ or $(10,5)$} \label{sec:B4_5}

Here the situation is almost the same as in
section~\ref{sec:b4-lines1}. In order to get $\widehat{X}$, we need to
blow up the image in $\P^3$ of the $G$-orbit of $L_{123}$ and the
$G$-orbit of $L_{234}$, then blow up the strict transform of the
$G$-orbit of $L_{12}$. 

The incidence data in Figure~\ref{fig:comb-B4} indicates
that $L_{12}$ does not contain any point in the $G$-orbit of
$L_{234}$. Indeed, each 2-plane in the $G$-orbit of $L_{12}$ contains
precisely 4 one-dimensional mirror intersections, two in the $G$-orbit
of $L_{123}$ and two in the $G$-orbit of $L_{124}$.

In other words, one gets the same formula as in
section~\ref{sec:b4-lines1} with $D$ replaced by $D_1$ and $D_2$, but
$D_2$ has no interaction with either $D_1$ or $E$. With the notation
\begin{eqnarray*}
  &\lambda=-4+4\nu_1+12\nu_2\\
  &\alpha_1=2+\delta_1-3\nu_1-6\nu_2, \qquad \alpha_2=2+\delta_2 - 6\nu_2\\
  &\beta=2(\delta_1-\nu_1),
\end{eqnarray*}
this gives
\begin{eqnarray*}
  &(K_Y+\Delta)^3 = \lambda^3 + 4\alpha_1^3 + 8\alpha_2^3 + 6\beta^3\cdot 2 
  - 3 \cdot 6\cdot \lambda\beta^2 
  - 3 \cdot 6\cdot 2\cdot \alpha_1\beta^2 ,
\end{eqnarray*}
For $p_1=6$, $p_2=6$, we take $\nu_j=1-\frac{2}{p_j}$ and
$\delta_1=2-\kappa_{L{123}}=0$, $\delta_2=2-\kappa_{L_{234}}=1-\frac{1}{3}$
(see the Table in Figure~\ref{fig:comb-B4}), which gives
$m_1=1$, $m_2=3$, we get
$$
  \chi^{orb}(\CHL(B_4,6,6))=-\frac{5}{432},
$$
and for $p_1=10$, $p_2=5$, $\delta_1=2-\kappa_{L{123}}=0$, $\delta_2=2-\kappa_{L_{234}}=1-\frac{1}{5}$
we get
$$
  \chi^{orb}(\CHL(B_4,10,5))=-\frac{13}{1500},
$$
as expected.


\subsubsection{The case $\CHL(B_4(2,8))$} \label{sec:B4_6}

This case is similar to the previous one. We now wish to compute
$$
\left(K_Y+\nu_1f_*M_1+\nu_2f_*M_2+\delta_1D_1+\delta_2D_2\right)^3,
$$
where $\nu_j=1-\frac{2}{p_j}$, $\delta_1=2-\kappa_{L_{123}}$, $\delta_2=2-\kappa_{L_{234}}$. Note that
\begin{eqnarray*}
  &K_{\widehat{X}}=K_Y\\
  &f^*f_*M_1=M_1+\frac{1}{2}E,\qquad f^*f_*M_2=M_2\\
  &f^*f_*D_1=D_1+\frac{1}{2}E,\qquad f^*f_*D_2=D_2+E.
\end{eqnarray*}
Indeed, each line in $\P^3$ below a component of $E$ contains three of
the points that get blown-up (one in the orbit of $L_{123}$, two in
the orbit of $L_{234}$), and it has a single transverse intersection
with a mirror in the first orbit of mirrors.

Using the blow-up map and the combinatorics of the arrangement, we have
\begin{eqnarray*}
  &\pi^*\pi_*M_1=M_1+3D_1\\
  &\pi^*\pi_*M_2=M_2+6D_1+6D_2+3E,
\end{eqnarray*}
and computations similar to the ones in the previous sections show
that $(K_Y+\Delta)^3$ is given by
$$
\left(\lambda\pi^*H+\alpha_1D_1+\alpha_2D_2+\beta E\right)^3
$$
where
\begin{eqnarray*}
  &\lambda=-4+4\nu_1+12\nu_2\\
  &\alpha_1=2+\delta_1-3\nu_1-6\nu_2,\qquad \alpha_2=2+\delta_2-6\nu_2\\
  &\beta=1+\frac{1}{2}\nu_1-3\nu_2+\frac{1}{2}\delta_1+\delta_2.
\end{eqnarray*}
Finally, developing the cube, we get
$$ 
\lambda^3 + 4\alpha_1^3 + 8 \alpha_2^3 + 16 \beta^4 \cdot 4 +
3\alpha_1\beta^2 D_1\cdot E^2 + 3\alpha_2\beta^2 D_2\cdot E^2 +
3\lambda\beta^2 \pi^*H\cdot E^2.
$$

Using the combinatorics and the above description for $E_j|_{E_j}$
(see equation~\eqref{eq:self}), we get
$$
\lambda^3 + 4\alpha_1^3 + 8 \alpha_2^3 + 16 \beta^4 \cdot 4 
- 3 \cdot 16 \cdot \alpha_1\beta^2 
- 3 \cdot 16\cdot 2\cdot \alpha_2\beta^2 
- 3 \cdot 16\cdot \lambda\beta^2.
$$
This gives 
$$
  \chi^{orb}(\CHL(B_4,2,8))=-\frac{11}{1024}=-\frac{11}{5120}\cdot 5,
$$ 
as it should since is has index 5 in the corresponding Deligne-Mostow
group (see Figure~\ref{tab:b4-dm}).


\subsubsection{The case $\CHL(B_4(4,8))$} \label{sec:B4_7}

This case is the most painful case to handle, but it simply combines
the difficulties we have encountered before. Here we blow up the
orbits of $L_{123}$ (4 copies), $L_{234}$ (8 copies), $L_{12}$ (6
copies) and $L_{23}$ (16 copies). Accordingly we have 4 exceptionals
$D_1,D_2,E_1,E_2$ in $\widehat{X}$, and still wish to compute
$$
\left(K_Y+\nu_1f_*M_1+\nu_2f_*M_2+\delta_1D_1+\delta_2D_2\right)^3,
$$
again with $\nu_j=1-\frac{2}{p_j}$, $\delta_1=2-\kappa_{L_{123}}$, $\delta_2=2-\kappa_{L_{234}}$.
Note that
\begin{eqnarray*}
  &K_{\widehat{X}}=f^*K_Y+E_1\\
  &f^*f_*M_1=M_1+\frac{1}{2}E_2,\qquad f^*f_*M_2=M_2+2E_1\\
  &f^*f_*D_1=D_1+2E_1+\frac{1}{2}E_2,\qquad f^*f_*D_2=D_2+E_2\\
  &K_{\widehat{X}}=\pi^*K_X+2D_1+2D_2+E_1+E_2\\
  &\pi^*\pi_*M_1=M_1+3D_1+2E_1,\qquad \pi^*\pi_*M_2=M_2+6D_1+6D_2+2E_1+3E_2
\end{eqnarray*}
The same computations as before now give
$$
\left(K_X+\Delta\right)^3=\left(\lambda\pi^*H+\alpha_1D_1+\alpha_2D_2+\beta_1E_1+\beta_2E_2\right)^3
$$
where
\begin{eqnarray*}
  &\lambda=-4+4\nu_1+12\nu_2\\
  &\alpha_1=2+\delta_1-3\nu_1-6\nu_2,\qquad
  \alpha_2=2+\delta_2-6\nu_2\\
  &\beta_1=2(\delta_1-\nu_2),\qquad
  \beta_2=1+\frac{1}{2}\nu_1-3\nu_2+\frac{1}{2}\delta_1+\delta_2.
\end{eqnarray*}
Inspecting the combinatorics of the arrangement and using
$E_1^{(j)}|_{E_1^{(j)}}=-l_1-2l_2$, $E_2^{(j)}|_{E_2^{(j)}}=-l_1-l_2$, we then get
$$
\lambda^3 + 4\alpha_1^3 + 8\alpha_2^3 + 16\cdot 4\cdot \beta_1^3 + 6\cdot 2\cdot \beta_2^3
- 3\cdot 6 \lambda\beta_1^2
- 3\cdot 16 \lambda\beta_2^2
- 3\cdot 6\cdot 2 \alpha_1\beta_1^2
- 3\cdot 16 \alpha_1\beta_2^2
- 3\cdot 16\cdot 2 \alpha_2\beta_2^2.
$$
This gives 
$$
  \chi^{orb}(\CHL(B_4,4,8))=-\frac{11}{1024},
$$ 
which is again the expected value.


\section{Presentations}\label{sec:presentations}

From the above results, one can easily obtain explicit presentations
for the CHL lattices. Indeed, recall that we denote $V=\C^{n+1}$,
$V^0\subset V$ the complement of the arrangement (given by the union
of the mirrors of reflections in $G$).  According to Theorem~7.1
in~\cite{chl}, a presentation for the linear holonomy group is given
by adjoining to a presentation of the braid group $\pi_1(G\backslash
V^0)$ some specific relations corresponding to the (irreducible)
strata in the arrangement. More specifically, for each irreducible
stratum $L$, consider the set of mirrors $\mathcal{H}_L$ that contain
$L$, and the braid group $G_L$ generated by the reflections in the
elements in $\mathcal{H}_L$, which has infinite cyclic center,
generated by an element $\alpha_L$. If $\rho:\pi_1(G\backslash
V^0)\rightarrow \Gamma$ denotes the holonomy representation, the CHL
relations correspond to imposing the order of $\rho(\alpha_L)$, given
by the integer that occurs in the Schwarz
condition~\eqref{eq:schwarz0}. In fact, among those relations, only
the ones where the mirror intersection $L$ of dimension or codimension
one are needed, since these are the such that the fixed point set of
the local holonomy group has fixed point set of codimension one (and
these are enough to present the orbifold fundamental group).

Presentations $\pi_1(\P(V^0/G))$ are given in~\cite{bessismichel}
(some of the results given there were conjectural at the time, but the
proof of their validity was given by Bessis
in~\cite{bessisannals}). It is easy to determine conjugacy classes of
loops corresponding to the conjugacy classes described in section 7.1
of~\cite{chl}, by determining the conjugacy classes of (irreducible)
mirror intersections in $G$, and then taking a generator of the center
of each stabilizer. 

The corresponding central elements are listed in Tables~1--5
in~\cite{brmaro}, for instance. One can also check their result by
using the explicit matrices described in~\cite{chl3d_1}. For example,
we list $(R_1R_2R_3)^3$ which generates the center of the braid group
generated by $R_1$, $R_2$ and $R_3$. Indeed, these generate a braid
group of type $G_{26}$, and a generator for the center is given in the
fifth column of Table~1 in~\cite{brmaro}.

We list the relevant central elements in
Table~\ref{fig:centralelements}; these give complex reflections in the
lattice, whose order is the integer occurring in the Schwarz condition
for $L$, and the relation is needed in the presentation only if
$\kappa_L>1$. For example, in the groups $\CHL(A_4,p)$,
$M=(R_1R_2R_3)^4$ is a complex reflection of order
$(\kappa_{L_{123}}-1)^{-1}=\frac{p}{p-4}$, and this relation is needed
in the presentation only for $p=5,6$ or $8$.
\begin{figure}
\begin{tabular}{|r|l|}
\hline
  $A_4$    & $(R_1R_2R_3)^4$\\
\hline
  $B_4$    & $(R_1R_2R_3)^3$, $(R_2R_3R_4)^4$\\
\hline
  $G_{28}$ & $(R_1R_2R_3)^3$, $(R_2R_3R_4)^3$\\
\hline
  $G_{29}$ & $(R_1R_2R_3)^4$, $(R_1R_2R_4)^3$, $(R_4R_3R_2)^8$, $(R_1R_2R_3^{-1}R_4R_3)^4$\\
\hline
  $G_{30}$ & $(R_1R_2R_3)^4$, $(R_2R_3R_4)^5$\\
\hline
  $G_{31}$ & $(R_5R_2R_1)^{\frac{4p}{p\wedge 3}}$, $(R_2R_3R_5)^4$\\
\hline
\end{tabular}
\caption{Complex reflections corresponding to central elements in
  $G_L$ for irreducible mirror intersection of dimension one.} \label{fig:centralelements}
\end{figure}

Collecting all this, we get the presentations in
Figure~\ref{fig:presentations} (p.~\pageref{fig:presentations}).
\begin{figure}
Presentation for $\CHL(A_4,p)$:\\
\begin{eqnarray*}
   \langle 
   \quad
   r_1,r_2,r_3,r_4  & | & r_1^p, (r_1r_2r_3)^{\frac{4p}{p-4}},\\
                    &   & \br_3(r_1,r_2), \br_3(r_2,r_3), \br_3(r_3,r_4), [r_1,r_3], [r_1,r_4], [r_2,r_4] 
   \quad 
   \rangle 
\end{eqnarray*}
{\centering \rule{8cm}{0.4pt}}\\
Presentation for $\CHL(B_4,p_1,p_2)$:\\
\begin{eqnarray*}
  \langle 
  \quad 
  r_1,r_2,r_3,r_4 & | & r_1^{p_1}, r_2^{p_2}, (r_1r_2r_3)^{\frac{3p_1p_2}{p_1p_2-2p_1-p_2}}, (r_2r_3r_4)^{\frac{4p_2}{p_2-4}},\\  
                  &   & \br_3(r_1,r_2), \br_3(r_2,r_3), \br_3(r_3,r_4), [r_1,r_3], [r_1,r_4], [r_2,r_4]
  \quad 
  \rangle 
\end{eqnarray*}\\
{\centering \rule{8cm}{0.4pt}}\\
\noindent Presentation for $\CHL(G_{28},p_1,p_2)$:\\
\begin{eqnarray*}
  \langle 
  \quad
  r_1,r_2,r_3,r_4 & | & r_1^{p_1}, r_3^{p_2}, (r_1r_2r_3)^{\frac{3p_1p_2}{p_1p_2-p_1-2p_2}}, (r_2r_3r_4)^{\frac{3p_1p_2}{p_1p_2-2p_1-p_2}},\\  
                  &   & \br_3(r_1,r_2), \br_4(r_2,r_3), \br_3(r_3,r_4), [r_1,r_3], [r_1,r_4], [r_2,r_4] 
  \quad 
  \rangle 
\end{eqnarray*}\\
{\centering \rule{8cm}{0.4pt}}\\
\noindent Presentation for $\CHL(G_{29},p)$:\\
\begin{eqnarray*}
  \langle 
  \quad
  r_1,r_2,r_3,r_4 & | & r_1^p, (r_1r_2r_4)^{\frac{3p}{p-3}}, (r_4r_3r_2)^{\frac{8p}{3p-8}},\\  
                  &   & \br_3(r_1,r_2), \br_4(r_2,r_3), \br_3(r_3,r_4), \br_4(r_2,r_4), \br_4(r_3,r_2r_4), [r_1,r_3], [r_1,r_4]
  \quad 
  \rangle 
\end{eqnarray*}\\
{\centering \rule{8cm}{0.4pt}}\\
\noindent Presentation for $\CHL(G_{30},p)$:\\
\begin{eqnarray*}
  \langle 
  \quad
  r_1,r_2,r_3,r_4 & | & r_1^p, (r_1r_2r_3)^{\frac{4p}{p-4}}, (r_2r_3r_4)^{\frac{10p}{4p-10}},\\  
                  &   & \br_3(r_1,r_2), \br_3(r_2,r_3), \br_5(r_3,r_4), [r_1,r_3], [r_1,r_4], [r_2,r_4]
  \quad 
  \rangle 
\end{eqnarray*}\\
{\centering \rule{8cm}{0.4pt}}\\
\noindent Presentation for $\CHL(G_{31},p)$:\\
\begin{eqnarray*}
  \langle 
  \quad
  r_1,r_2,r_3,r_4,r_5 & | & r_1^p,                     
                            (r_5r_2r_1)^{\frac{8p^2}{(p\wedge 3)(4p-10)}}, 
                            (r_2r_3r_5)^{\frac{4p}{p-4}}, 
                            r_1r_5r_4=r_5r_4r_1=r_4r_1r_5\\  
                      &   &  \br_3(r_1,r_2), \br_3(r_2,r_5), \br_3(r_5,r_3), \br_3(r_3,r_4), 
                             [r_1,r_3], [r_2,r_3], [r_2,r_4]
  \quad 
  \rangle 
\end{eqnarray*}
\caption{Presentations for Couwenberg-Heckman-Looijenga lattices in $PU(3,1)$}\label{fig:presentations}
\end{figure}

\section{Combinatorial data}\label{sec:comb}

\begin{figure}
  \includegraphics{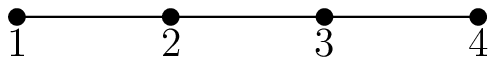}
  \ \\[1cm]

  \begin{tabular}{|c|c|c|c|c|c|}
    \hline
    $G$           & $|G|$ & $|Z(G)|$ & Mirror orbit & $|$orbit$|$ & Weight \\
    \hline
    $W(A_4),S_5$     &  120  &     1    &  $L_1$       &    10       &   $1-\frac{2}{p}$\\
    \hline
  \end{tabular}
  \ \\[1cm]
  
  \begin{tabular}{|c|c|c|}
    \hline
    Finite                        &   Parabolic   &   Hyperbolic\\
    \hline
    $p=2 (W(A_4))$, $3 (G_{32})$       &               &   $p=4,5,6,8$\\
    \hline
  \end{tabular}
  \ \\[1cm]
  
  \begin{tabular}{|c|c|c|c|c|c|c|}
    \hline
    $L$             & \#(orbit)            & \#(mirrors) & $|Z(G_L)|$       & $\kappa_L$                     &  Incident vertices\\
    \hline
    $L_{12}$        & 10                   &   3         &     1            & $\frac{3}{2}(1-\frac{2}{p})$   &  $2\times L_{123},L_{134}$\\
    $L_{13}$        & 15                   &   2         &                  & (reducible)                    &  $L_{123},2\times L_{134}$\\
    \hline
  \end{tabular}
  
  \begin{tabular}{|c|c|c|c|c|c|c|}
    \hline
    $p$                  & 2  &   3             &   4                  &   5              &      6               &        8  \\
    \hline
    $\kappa_{L_{12}}$   & 0  & $1-\frac{1}{2}$ & $1-\frac{1}{4}$      & $1-\frac{1}{10}$ & $1$                  &  $1+\frac{1}{8}$      \\
    \hline
  \end{tabular}
  \ \\[1cm]
  
  \begin{tabular}{|c|c|c|c|c|c|}
    \hline
    $L$           & \#(orbit)             & \#(mirrors) & $|Z(G_L)|$  & $\kappa_L$                            &  Incident lines\\
    \hline
    $L_{123}$     & 5                     &   6         &     1       & $2(1-\frac{2}{p})$    & $4\times L_{12},3\times L_{13}$\\
    $L_{134}$     & 10                    &   4         & (reducible) & (reducible)           & $1\times L_{12},3\times L_{13}$\\
    \hline
  \end{tabular}
  
  \begin{tabular}{|c|c|c|c|c|c|c|}
    \hline
    $p$                & 2    &    3            & 4     & 5               & 6               & 8     \\
    \hline
    $\kappa_{L_{123}}$ & 0    & $1-\frac{1}{3}$ & $1$   & $1+\frac{1}{5}$ & $1+\frac{1}{3}$ & $1+\frac{1}{2}$\\
    \hline
  \end{tabular}
\caption{Combinatorial data for $A_4$}\label{fig:comb-A4}
\end{figure}

\begin{figure}
  \includegraphics{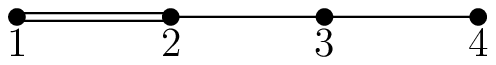}
  \ \\[1cm]
  
  \begin{tabular}{|c|c|c|c|c|c|}
    \hline
    $G$             & $|G|$ & $|Z(G)|$ & Mirror orbit    & $|$orbit$|$ & Weight \\
    \hline
    $W(B_4),G(2,1,4)$  &  384  &    2     &  $L_1$           &    4        &   $1-\frac{2}{p_1}$\\
                       &       &          &  $L_2$           &    12       &   $1-\frac{2}{p_2}$\\
    \hline
  \end{tabular}
  \ \\[1cm]
  
  {\footnotesize
    \begin{tabular}{|c|c|c|}
      \hline
      Finite                        &   Parabolic         &   Hyperbolic\\
      \hline
      $(p_1,p_2)=(n,2)(G(n,1,4))$   &   $(p_1,p_2)=(2,3)$ &   $(p_1,p_2)=(2,4),(2,5),(2,6),(2,8),(3,3),(3,4),(3,6)$\\
      &                     &   $(4,3),(4,4),(4,8),(6,3),(6,4),(6,6),(10,5),(12,3)$\\
      \hline  
  \end{tabular}}
  \ \\[1cm]
  \begin{rk}
    The group derived from $B_4$ and orders $(p_1,p_2)=(5,5)$ is the
    Deligne-Mostow group $(3,3,3,3,3,5)/10$, so it is a lattice; however
    it does not satifsy the Schwarz condition in~\cite{chl}, since in
    that case $\kappa_{L_{123}}-1=4/5$, but $|Z(G_{L_{123}})|=2$ only
    allows numerator 1 or 2, not 4. This group can also be described as
    $\CHL(A_4,5)$, where the Schwarz condition does hold.
  \end{rk}
  
  \begin{tabular}{|c|c|c|c|c|c|}
    \hline
    $L$             & \#(orbit)   & \#(mirrors) & $|Z(G_L)|$       & $\kappa_L$                            &  Incident vertices\\
    \hline
    $L_{12}$        & 6           &   2+2       &     2            & $(1-\frac{2}{p_1})+(1-\frac{2}{p_2})$ & $2\times L_{123},2\times L_{124}$\\
    $L_{14}$        & 24          &   1+1       &  (reducible)     & (reducible)                           & $L_{123},L_{124},2\times L_{134}$\\
    $L_{23}$        & 16          &   0+3       &     1            & $\frac{3}{2}(1-\frac{2}{p_2})$        & $L_{123},L_{134},2\times L_{234}$\\
    $L_{24}$        & 12          &   0+2       &  (reducible)     & (reducible)                           & $2\times L_{124},2\times L_{234}$\\
    \hline
  \end{tabular}
  \begin{tabular}{|c|c|c|c|c|c|c|c|}
    \hline
    $(p_1,p_2)$       & (2,3)           & (2,4)           & (2,5)            & (2,6)           & (2,8)        \\
    \hline
    $\kappa_{L_{12}}$ & $1-\frac{2}{3}$ & $1-\frac{1}{2}$ & $1-\frac{2}{5}$  & $1-\frac{1}{3}$ & $1-\frac{1}{4}$           \\
    $\kappa_{L_{23}}$ & $1-\frac{1}{2}$ & $1-\frac{1}{4}$ & $1-\frac{1}{10}$ & $1$             & $1+\frac{1}{8}$ \\
    \hline
  \end{tabular}
  \begin{tabular}{|c|c|c|c|c|c|c|c|}
    \hline
    $(p_1,p_2)$       & (3,3)           & (3,4)           & (3,6)  & (4,3)             & (4,4)           &    (4,8) \\
    \hline
    $\kappa_{L_{12}}$ & $1-\frac{1}{3}$ & $1-\frac{1}{6}$ & $1$    & $1-\frac{1}{6}$   & $1$             & $1+\frac{1}{4}$\\
    $\kappa_{L_{23}}$ & $1-\frac{1}{2}$ & $1-\frac{1}{4}$ & $1$    & $1-\frac{1}{2}$   & $1-\frac{1}{4}$ & $1+\frac{1}{8}$\\
    \hline
  \end{tabular}
  \begin{tabular}{|c|c|c|c|c|c|c|c|}
    \hline
    $(p_1,p_2)$       & (6,3)           & (6,4)           & (6,6)           & (10,5)           & (12,3)      \\
    \hline
    $\kappa_{L_{12}}$ & $1$             & $1+\frac{1}{6}$ & $1+\frac{1}{3}$ & $1+\frac{2}{5}$  & $1+\frac{1}{6}$  \\
    $\kappa_{L_{23}}$ & $1-\frac{1}{2}$ & $1-\frac{1}{4}$ & $1$             & $1-\frac{1}{10}$ & $1-\frac{1}{2}$  \\
    \hline
  \end{tabular}
  \ \\[1cm]
  
  \begin{tabular}{|c|c|c|c|c|c|}
    \hline
    $L$        & \#(orbit) & \#(mirrors) &    $|Z(G_L)|$   &   $\kappa_L$                          &  Incident lines\\
    \hline
    $L_{123}$  & 4       & 3+6         &        2        & $1-\frac{2}{p_1}+2(1-\frac{2}{p_2})$  &  $3\times L_{12},4\times L_{23},6\times L_{14}$\\
    $L_{124}$  & 12      & 2+3         & (reducible)     & (reducible)                           &  $1\times L_{12},2\times L_{14},2\times L_{24}$\\
    $L_{134}$  & 16      & 1+3         & (reducible)     & (reducible)                           &  $3\times L_{14},1\times L_{23}$\\
    $L_{234}$  & 8       & 0+6         &        1        & $2(1-\frac{2}{p_2})$                  &  $4\times L_{23},3\times L_{24}$\\
    \hline
  \end{tabular}
  \begin{tabular}{|c|c|c|c|c|c|c|c|}
    \hline
    $(p_1,p_2)$        & (2,3)           & (2,4)  & (2,5)            & (2,6)           & (2,8)        \\
    \hline
    $\kappa_{L_{123}}$ & $1-\frac{1}{3}$ & $1$    & $1+\frac{1}{5}$  & $1+\frac{1}{3}$ & $1+\frac{1}{2}$           \\
    $\kappa_{L_{234}}$ & $1-\frac{1}{3}$ & $1$    & $1+\frac{1}{5}$  & $1+\frac{1}{3}$ & $1+\frac{1}{2}$ \\
    \hline
  \end{tabular}
  \begin{tabular}{|c|c|c|c|c|c|c|c|}
    \hline
    $(p_1,p_2)$        & (3,3)           & (3,4)           & (3,6)           & (4,3)             & (4,4)           &    (4,8) \\
    \hline
    $\kappa_{L_{123}}$ & $1$             & $1+\frac{1}{3}$ & $1+\frac{2}{3}$ & $1+\frac{1}{6}$   & $1+\frac{1}{2}$ & $1+\frac{1}{1}$\\
    $\kappa_{L_{234}}$ & $1-\frac{1}{3}$ & $1$             & $1+\frac{1}{3}$ & $1-\frac{1}{3}$   & $1$             & $1+\frac{1}{2}$\\
    \hline
  \end{tabular}
  \begin{tabular}{|c|c|c|c|c|c|c|c|}
    \hline
    $(p_1,p_2)$        & (6,3)           & (6,4)           & (6,6)           & (10,5)           & (12,3)      \\
    \hline
    $\kappa_{L_{123}}$ & $1+\frac{1}{3}$ & $1+\frac{2}{3}$ & $1+\frac{1}{1}$ & $1+\frac{1}{1}$  & $1+\frac{1}{2}$  \\
    $\kappa_{L_{234}}$ & $1-\frac{1}{3}$ & $1$             & $1+\frac{1}{3}$ & $1+\frac{1}{5}$  & $1-\frac{1}{3}$  \\
    \hline
  \end{tabular}
\caption{Combinatorial data for $B_4$}\label{fig:comb-B4}
\end{figure}

\begin{figure}
    \includegraphics{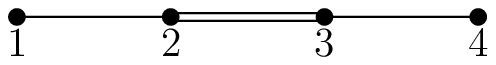}
    \ \\[1cm]
    
    \begin{tabular}{|c|c|c|c|c|c|}
      \hline
      $G$           & $|G|$ & $|Z(G)|$ & Mirror orbit & $|$orbit$|$ & Weight \\
      \hline
      $G_{28}, W(F_4)$ &  1152 &     2    &  $L_1$    &    12     &   $1-\frac{2}{p_1}$\\
                       &       &          &  $L_3$    &    12     &   $1-\frac{2}{p_2}$\\
      \hline
    \end{tabular}
    \ \\[1cm]
    
    \begin{tabular}{|c|c|c|}
      \hline
      Finite                        &   Parabolic         &   Hyperbolic\\
      \hline
      $(p_1,p_2)=(2,2)(G_{28})$     &   $(p_1,p_2)=(2,3)$ &   $(p_1,p_2)=(2,4),(2,5),(2,6),(2,8),(2,12)$\\
      &                     &   $(3,3),(3,4),(3,6),(3,12),(4,4),(6,6)$\\
      \hline
    \end{tabular}
    \ \\[1cm]
    
    \begin{tabular}{|c|c|c|c|c|c|}
      \hline
      $L$       & \#(orbit)   & \#(mirrors) &     $|Z(G_L)|$              &   $\kappa_L$                              & Incident vertices\\
      \hline
      $L_{12}$  & 16        & 3+0         &                   1         & $3(\frac{1}{2}-\frac{1}{p_1})$            & $3\times L_{123}, 3\times L_{124}$\\
      $L_{14}$  & 72        & 1+1         & (reducible)                 & (reducible)                               & $L_{123}, L_{234}, 2\times L_{134}, 2\times L_{124}$\\
      $L_{23}$  & 18        & 2+2         &                   2         & $2(1-\frac{1}{p_1}-\frac{1}{p_2})$        & $2\times L_{123},2\times L_{234}$\\
      $L_{34}$  & 16        & 0+3         &                   1         & $3(\frac{1}{2}-\frac{1}{p_2})$            & $3\times L_{234},3\times L_{134}$\\
      \hline
    \end{tabular}
    
    \begin{tabular}{|c|c|c|c|c|c|c|c|}
      \hline
      $(p_1,p_2)$       & (2,2) & (2,3)       & (2,4)       & (2,5)       & (2,6)       & (2,8)         & (2,12)\\
      \hline
      $\kappa_{L_{12}}$ & 0     &  0          &   0         &   0         &   0         &   0           &   0\\
      $\kappa_{L_{23}}$ & 0     & $\frac{1}{3}$ & $\frac{1}{2}$ & $\frac{3}{5}$ & $\frac{2}{3}$ & $\frac{3}{4}$   & $\frac{5}{6}$\\
      $\kappa_{L_{34}}$ & 0     & $\frac{1}{2}$ & $\frac{3}{4}$ & $\frac{9}{10}$&  1          & $1+\frac{1}{8}$ & $1+\frac{1}{4}$\\
      \hline
    \end{tabular}
    \begin{tabular}{|c|c|c|c|c|c|c|}
      \hline
      $(p_1,p_2)$       & (3,3)       & (3,4)           & (3,6)         & (3,12)          & (4,4)         & (6,6) \\
      \hline
      $\kappa_{L_{12}}$ & $\frac{1}{2}$ & $\frac{1}{2}$ & $\frac{1}{2}$ & $\frac{1}{2}$   & $\frac{3}{4}$ & 1\\
      $\kappa_{L_{23}}$ & $\frac{2}{3}$ & $\frac{5}{6}$ & 1             & $1+\frac{1}{6}$ & 1             & $1+\frac{1}{3}$\\
      $\kappa_{L_{34}}$ & $\frac{1}{2}$ & $\frac{3}{4}$ & 1             & $1+\frac{1}{4}$ & $\frac{3}{4}$ & 1\\
      \hline
    \end{tabular}
    \ \\[1cm]
    
    \begin{tabular}{|c|c|c|c|c|c|}
      \hline
      $L$        & \#(orbit) & \#(mirrors) &    $|Z(G_L)|$   &   $\kappa_L$                             & Incident lines \\
      \hline
      $L_{123}$  & 12      & 6+3         &        2        & $2(1-\frac{2}{p_1})+(1-\frac{2}{p_2})$   & $4\times L_{12},6\times L_{14},3\times L_{23}$\\
      $L_{234}$  & 12      & 3+6         &        2        & $(1-\frac{2}{p_1})+2(1-\frac{2}{p_2})$   & $6\times L_{14},3\times L_{23},4\times L_{34}$\\
      $L_{134}$  & 48      & 1+3         & (reducible)     & (reducible)                              & $3\times L_{14},1\times L_{34}$\\
      $L_{124}$  & 48      & 3+1         & (reducible)     & (reducible)                              & $1\times L_{12},3\times L_{14}$\\
      \hline
    \end{tabular}
    \begin{tabular}{|c|c|c|c|c|c|c|c|}
      \hline
      $(p_1,p_2)$        & (2,2) & (2,3)         & (2,4)         & (2,5)           & (2,6)           & (2,8)            & (2,12)\\
      \hline
      $\kappa_{L_{123}}$ & 0     & $\frac{1}{3}$ & $\frac{1}{2}$ & $\frac{3}{5}$   & $\frac{2}{3}$   & $\frac{3}{4}$    & $\frac{5}{6}$\\
      $\kappa_{L_{234}}$ & 0     & $\frac{2}{3}$ &    1          & $1+\frac{1}{5}$ & $1+\frac{1}{3}$ & $1+\frac{1}{2}$  & $1+\frac{2}{3}$\\
      \hline
    \end{tabular}
    \begin{tabular}{|c|c|c|c|c|c|c|}
      \hline
      $(p_1,p_2)$        & (3,3) & (3,4)           & (3,6)           & (3,12)          & (4,4)           & (6,6) \\
      \hline
      $\kappa_{L_{123}}$ &   1   & $1+\frac{1}{6}$ & $1+\frac{1}{3}$ & $1+\frac{1}{2}$ & $1+\frac{1}{2}$ & $1+1$\\
      $\kappa_{L_{234}}$ &   1   & $1+\frac{1}{3}$ & $1+\frac{2}{3}$ & $1+1$           & $1+\frac{1}{2}$ & $1+1$\\
      \hline
    \end{tabular}
\caption{Combinatorial data for $G_{28}$.}\label{fig:comb-g28}
\end{figure}

\begin{figure}
  \includegraphics{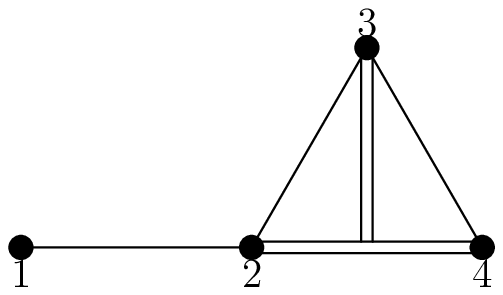}
  \ \\[1cm]
  
  \begin{tabular}{|c|c|c|c|}
    \hline
    $G$           & $|G|$ & $|Z(G)|$ & \#(mirrors)\\
    \hline
    $G_{29}$      & 7680  &     4    &    40 \\
    \hline
  \end{tabular}
  \ \\[1cm]
  
  \begin{tabular}{|c|c|c|}
    \hline
    Finite                        &   Parabolic   &   Hyperbolic\\
    \hline
    $p=2(G_{29})$                 &               &   $p=3,4$\\
    \hline
  \end{tabular}
  \ \\[1cm]
  
  \begin{tabular}{|c|c|c|c|c|c|}
    \hline
    $L$             & \#(orbit)   & \#(mirrors) & $|Z(G_L)|$       & $\kappa_L$                    &  Incident vertices \\
    \hline
    $L_{12}$        & 160       &   3         &     1            & $\frac{3}{2}(1-\frac{2}{p})$  & $2\times L_{123},2\times L_{12343},L_{124},L_{134},2\times L_{234}$\\
    $L_{13}$        & 120       &   2         &                  & (reducible)                   & $2\times L_{123},2\times L_{12343},2\times L_{124},4\times L_{134}$\\
    $L_{24}$        & 30        &   4         &     2            & $2(1-\frac{2}{p})$            & $4\times L_{12343},2\times L_{234}$\\
    \hline
  \end{tabular}
  
  \begin{tabular}{|c|c|c|c|}
    \hline
    $p$                  & 2  &   3           &   4  \\
    \hline
    $\kappa_{L_{12}}$   & 0  & $\frac{1}{2}$ & $\frac{3}{4}$             \\
    $\kappa_{L_{24}}$   & 0  & $\frac{2}{3}$ & $1$ \\
    \hline
  \end{tabular}
  \ \\[1cm]
  
  \begin{tabular}{|c|c|c|c|c|c|}
    \hline
    $L$           & \#(orbit)  & \#(mirrors) & $|Z(G_L)|$  & $\kappa_L$         & Incident lines\\
    \hline
    $L_{123}$     & 80      &   6         &     1       & $2(1-\frac{2}{p})$  & $4\times L_{12},3\times L_{13}$\\
    $L_{12343}$   & 80      &   6         &     1       & $2(1-\frac{2}{p})$  & $4\times L_{12},3\times L_{13}$\\
    $L_{124}$     & 40      &   9         &     2       & $3(1-\frac{2}{p})$  & $4\times L_{12},6\times L_{13},3\times L_{24}$\\
    $L_{134}$     & 160     &   4         & (reducible) & (reducible)         & $1\times L_{12},3\times L_{13}$\\
    $L_{234}$     & 20      &   12        &     1       & $4(1-\frac{2}{p})$  & $16\times L_{12},3\times L_{24}$\\
    \hline
  \end{tabular}
  
  \begin{tabular}{|c|c|c|}
    \hline
    $p$                & 3               & 4         \\
    \hline
    $\kappa_{L_{123}}$ & $\frac{2}{3}$   & $1$ \\
    $\kappa_{L_{124}}$ & $1$             & $1+\frac{1}{2}$ \\
    $\kappa_{L_{234}}$ & $1+\frac{1}{3}$ & $1+1$ \\
    \hline
  \end{tabular}
\caption{Combinatorial data for $G_{29}$}\label{fig:comb-g29}
\end{figure}

\begin{figure}
  \includegraphics{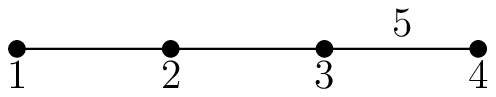}
  \ \\[1cm]

  \begin{tabular}{|c|c|c|c|}
    \hline
    $G$               & $|G|$ & $|Z(G)|$ & \#(mirrors)\\
    \hline
    $G_{30}$, $W(H_4)$   & 14400 &     2    &    60 \\
    \hline
  \end{tabular}
  \ \\[1cm]
  
  \begin{tabular}{|c|c|c|}
    \hline
    Finite                        &   Parabolic   &   Hyperbolic\\
    \hline
    $p=2(G_{30})$                 &               &   $p=3,5$\\
    \hline
  \end{tabular}
  \ \\[1cm]
  
  \begin{tabular}{|c|c|c|c|c|c|}
    \hline
    $L$             & \#(orbit)   & \#(mirrors) & $|Z(G_L)|$       & $\kappa_L$                    &  Incident vertices\\
    \hline
    $L_{12}$        & 200       &   3         &     1            & $\frac{3}{2}(1-\frac{2}{p})$  & $6\times L_{123},3\times L_{124},3\times L_{234}$\\
    $L_{13}$        & 450       &   2         &  (reducible)     & (reducible)                   & $2\times L_{123},4\times L_{124},4\times L_{134},2\times L_{234}$\\
    $L_{34}$        & 72        &   5         &     1            & $\frac{5}{2}(1-\frac{2}{p})$  & $5\times L_{134},5\times L_{234}$\\
    \hline
  \end{tabular}
  
  \begin{tabular}{|c|c|c|c|}
    \hline
    $p$                  & 2  &   3           &   5  \\
    \hline
    $\kappa_{L_{12}}$   & 0  & $\frac{1}{2}$ & $\frac{9}{10}$ \\
    $\kappa_{L_{34}}$   & 0  & $\frac{5}{6}$ & $1+\frac{1}{2}$ \\
    \hline
  \end{tabular}
  \ \\[1cm]
  
  \begin{tabular}{|c|c|c|c|c|c|}
    \hline
    $L$           & \#(orbit)  & \#(mirrors) & $|Z(G_L)|$  & $\kappa_L$         &  Incident lines\\
    \hline
    $L_{123}$     & 300     &   6         &     1       & $2(1-\frac{2}{p})$  & $4\times L_{12}, 3\times L_{13}$\\
    $L_{124}$     & 600     &   4         & (reducible) & (reducible)         & $1\times L_{12}, 3\times L_{13}$\\
    $L_{134}$     & 360     &   6         & (reducible) & (reducible)         & $5\times L_{13}, 1\times L_{34}$\\
    $L_{234}$     & 60      &   15        &     2       & $5(1-\frac{2}{p})$  & $10\times L_{12}, 15\times L_{13}, 6\times L_{34}$\\
    \hline
  \end{tabular}
  
  \begin{tabular}{|c|c|c|c|c|c|c|c|}
    \hline
    $p$                & 3               & 5         \\
    \hline
    $\kappa_{L_{123}}$ & $\frac{2}{3}$   & $1+\frac{1}{5}$ \\
    $\kappa_{L_{234}}$ & $1+\frac{2}{3}$ & $1+\frac{2}{1}$ \\
    \hline
  \end{tabular}
\caption{Combinatorial data for $G_{30}$}\label{fig:comb-g30}
\end{figure}

\begin{figure}
  \includegraphics{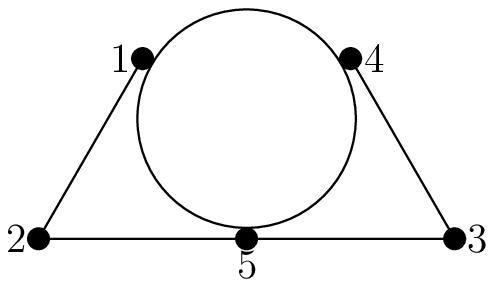}
  \ \\[1cm]
  
  \begin{tabular}{|c|c|c|c|}
    \hline
    $G$           & $|G|$ & $|Z(G)|$ & \#(mirrors)\\
    \hline
    $G_{31}$      & 46080 &     4    &    60 \\
    \hline
  \end{tabular}
  \ \\[1cm]
  
  \begin{tabular}{|c|c|c|}
    \hline
    Finite                        &   Parabolic   &   Hyperbolic\\
    \hline
    $p=2(G_{31})$                 &               &   $p=3,5$\\
    \hline
  \end{tabular}
  \ \\[1cm]
  
  \begin{tabular}{|c|c|c|c|c|c|}
    \hline
    $L$                  & \#(orbit)   & \#(mirrors) & $|Z(G_L)|$       & $\kappa_L$                     &  Incident vertices\\
    \hline
    $L_{12}$        & 320       &   3         &     1            & $\frac{3}{2}(1-\frac{2}{p})$   & $3\times L_{123},3\times L_{125},6\times L_{235}$\\
    $L_{13}$        & 360       &   2         &  (reducible)     & (reducible)                    & $8\times L_{123},2\times L_{125},4\times L_{235}$\\
    $L_{14}$        & 30        &   6         &     4            & $3(1-\frac{2}{p})$             & $6\times L_{125}$\\
    \hline
  \end{tabular}
  
  \begin{tabular}{|c|c|c|c|}
    \hline
    $p$                  & 2  &   3           &   5  \\
    \hline
    $\kappa_{L_{12}}$   & 0  & $\frac{1}{2}$ & $\frac{9}{10}$ \\
    $\kappa_{L_{14}}$   & 0  & $1$           & $1+\frac{4}{5}$ \\
    \hline
  \end{tabular}
  \ \\[1cm]
  
  \begin{tabular}{|c|c|c|c|c|c|}
    \hline
    $L$           & \#(orbit)  & \#(mirrors) & $|Z(G_L)|$  & $\kappa_L$         &  Incident vertices\\
    \hline
    $L_{123}$     & 960     &   4         & (reducible) & (reducible)         & $1\times L_{12}, 3\times L_{13}$\\
    $L_{125}$     & 60      &   15        &       2     & $5(1-\frac{2}{p})$  & $16\times L_{12}, 12\times L_{13}, 3\times L_{14}$\\
    $L_{235}$     & 480     &   6         &       1     & $2(1-\frac{2}{p})$  & $4\times L_{12}, 3\times L_{13}$\\
    \hline
  \end{tabular}
  \begin{tabular}{|c|c|c|c|c|c|c|c|}
    \hline
    $p$                & 3               & 5         \\
    \hline
    $\kappa_{L_{125}}$ & $1+\frac{2}{3}$ & $1+\frac{2}{1}$ \\
    $\kappa_{L_{235}}$ & $\frac{2}{3}$   & $1+\frac{1}{5}$ \\
    \hline
  \end{tabular}
\caption{Combinatorial data for $G_{31}$}\label{fig:comb-g31}
\end{figure}

In Figures~\ref{fig:comb-A4} through~\ref{fig:comb-g31}
(pp.~\pageref{fig:comb-A4}-\pageref{fig:comb-g31}), we list combinatorial
data that allow us to check the Schwarz conditions (see section~4
of~\cite{chl}) and to compute volumes (see
section~\ref{sec:volumes}). 

For the group $G_{28}$, there are two
orbits of mirrors, which can be assigned independent
weights. Accordingly, we give the number of mirrors containing a given
$L$ in the form $j+k$, where $j$ (resp. $k$) is the number of mirrors
from the first (resp. second) orbit.

For each group orbit of irreducible mirror intersections (see p.~88
of~\cite{chl}), we list the corresponding weight $\kappa_L$, which is
the ratio
\begin{equation}\label{eq:schwarz}
  \kappa_L=\frac{\sum_{H\in H_L} \kappa_H}{{\rm codim} L},
\end{equation}
where $H_L$ is the set of hyperplanes in the mirror arrangement
that contain $L$.

We also list the order of the center $Z(G_L)$ of the Schwarz symmetry
group $G_L$. Recall that $G_L$ is obtained as the fixed point
stabilizer of $L$, and it is a reflection group (generated by the
reflections in $G$ whose mirror contains $L$).

The Schwarz condition amounts to requiring that, for every irreducible
mirror intersection $L$ such that $\kappa_L>1$,
$$
\kappa_L-1=\frac{|Z(G_L)|}{n_L}
$$
for some integer $n_L\geq 2$.

Since the condition applies only to irreducible mirror intersections,
when $L$ is not irreducible, we do not compute any weight, and simply
write ``(reducible)'' in the corresponding spot in the table.

In order to describe strata in the arrangement, we label them with an
subscript that indicates the mirrors of reflections that define a given
intersection using the numbering of the reflection generators. For
instance, $L_j$ denotes the mirror of the $j$-th reflection $R_j$,
$L_{jk}$ denotes the intersection of the mirrors of the reflections
$R_j$ and $R_k$, $L_{ijk}$ denotes the intersection of the three
mirrors of $R_i$, $R_j$ and $R_k$, etc. We extend this notation
slightly to include conjugates of the generators, for instance
$L_{12343}$ denotes the intersection of the mirrors of $R_1$, $R_2$
and $R_3R_4R_3$.

When computing volumes, we will need some data on incidence relations
between mirror intersections of various dimensions; what we need is
listed in the columns with header ``Incident vertices'' or ``Incident
lines''. Recall that vertices (resp. lines) in $\P V$
actually correspond to lines (resp. 2-planes) in $V$.

When we write ``$2\times L_{123},L_{134}$'' in the column for incident
vertices to $L_{12}$ (see the table in Figure~\ref{fig:comb-A4} for
the $A_4$ arrangement), we mean that $L_{12}$ contains three
1-dimensional mirror intersections, and among those three, two that
are in the $G$-orbit of $L_{123}$ and one is in the orbit of
$L_{134}$. We only use this notation provided the $G$-orbits of
$L_{123}$ and $L_{134}$ are disjoint.

\section{Volumes and rough commensurability invariants}

In tables~\ref{tab:b4-dm}~(p.~\pageref{tab:b4-dm})
and~\ref{fig:rough_invariants}~(p.~\pageref{fig:rough_invariants}), we
collect rough commensurability invariants (cocompactness,
arithmeticity, adjoint trace fields) and orbifold Euler characteristic
of CHL lattices. For groups known to be commensurable with
Deligne-Mostow lattices, we give the exponents of the relevant
hypergeometric functions, and the index in the corresponding maximal
Deligne-Mostow lattice.

\begin{table}[htbp]
\begin{tabular}{|c|c|c|c|c|c|c|c|}
\hline
  ST group & Order(s) & DM group           &  Index    &  C/NC  &   A/NA   &   $\tr Ad\Gamma$    & $\chi^{orb}$\\
\hline
  $A_4$    & 4        & $(1,1,1,1,1,3)/4$  &  1        &   NC   &    A     &      $\Q$           & $-1/1920$\\
           & 5        & $(3,3,3,3,3,5)/10$ &  1        &   C    &    A     &      $\Q(\sqrt{5})$ & $-{17}/{6000}$\\
           & 6        & $(1,1,1,1,1,1)/3$  &  6        &   NC   &    A     &    $\Q$             & $-1/270$\\
           & 8        & $(1,3,3,3,3,3)/8$  &  1        &   C    &    A     &    $\Q(\sqrt{2})$   & $-11/5120$\\
\hline
  $B_4$    & $(2,4)$    & $(1,1,1,1,1,3)/4$  &  5      &   NC   &    A     & $\Q$                & $-{1}/{384}$\\
           & $(2,5)$    & $(2,3,3,3,3,6)/10$ &  1      &   C    &    A     & $\Q(\sqrt{5})$      & $-{13}/{1500}$\\
           & $(2,6)$    & $(1,2,2,2,2,3)/6$  &  1      &   NC   &    A     & $\Q$                & $-{5}/{432}$\\
           & $(2,8)$    & $(1,3,3,3,3,3)/8$  &  5      &   C    &    A     & $\Q(\sqrt{2})$      & $-{11}/{1024}$ \\
           & $(3,3)$    & $(1,1,1,1,3,5)/6$  &  1      &   NC   &    A     & $\Q$                & $-{1}/{1296}$\\
           & $(3,4)$    & $(3,3,3,3,5,7)/12$ &  1      &   NC   &    NA    & $\Q(\sqrt{3})$      & $-{31}/{3456}$\\
           & $(3,6)$    & $(1,1,1,1,1,1)/3$  &  30     &   NC   &    A     & $\Q$                & $-{1}/{54}$\\
           & $(4,3)$    & $(2,2,2,2,7,9)/12$ &  1      &   C    &    A     & $\Q(\sqrt{3})$      & $-{23}/{10368}$\\
           & $(4,4)$    & $(1,1,1,1,2,2)/4$  &  2      &   NC   &    A     & $\Q$                & $-{1}/{96}$\\
           & $(4,8)$    & $(1,3,3,3,3,3)/8$  &  5      &   C    &    A     & $\Q(\sqrt{2})$      & $-{11}/{1024}$\\
           & $(6,3)$    & $(1,1,1,1,4,4)/6$  &  2      &   NC   &    A     & $\Q$                & $-1/324$\\
           & $(6,4)$    & $(3,3,3,3,5,7)/12$ &  1      &   NC   &    NA    & $\Q(\sqrt{3})$      & $-{31}/{3456}$\\
           & $(6,6)$    & $(1,2,2,2,2,3)/6$  &  1      &   NC   &    A     & $\Q$                & $-{5}/{432}$\\
           & $(10,5)$   & $(2,3,3,3,3,6)/10$ &  1      &   C    &    A     & $\Q(\sqrt{5})$      & $-{13}/{1500}$\\
           & $(12,3)$   & $(2,2,2,2,7,9)/12$ &  1      &   C    &    A     & $\Q(\sqrt{3})$      & $-{23}/{10368}$\\
\hline
\end{tabular}
\caption{CHL groups for $A_4$ and $B_4$ are subgroups of specific
  Deligne-Mostow lattices.}\label{tab:b4-dm}
\end{table}

\begin{table}[htbp]
\begin{tabular}{|c|c|c|c|c|c|c|}
  \hline
  ST Group           & $p$, $(p_1,p_2)$    & C/NC      &     A/NA   & Adjoint trace field   & Euler char.    \\
  \hline
  $G_{28}$           &   $(2,4)$           &     NC    &      A     &  $\Q$                 &  $-1/1152$     \\
                     &   $(2,5)$           &     C     &      A     &  $\Q(\sqrt{5})$       &  $-13/4500$    \\
                     &   $(2,6)$           &     NC    &      A     &  $\Q$                 &  $-5/1296$     \\
                     &   $(2,8)$           &     C     &      A     &  $\Q(\sqrt{2})$       &  $-11/3072$    \\
                     &   $(2,12)$          &     C     &      A     &  $\Q(\sqrt{3})$       &  $-23/10368$   \\
                     &   $(3,3)$           &     NC    &      A     &  $\Q$                 &  $-1/144$      \\
                     &   $(3,4)$           &     C     &      A     &  $\Q(\sqrt{3})$       &  $-23/1152$    \\
                     &   $(3,6)$           &     NC    &      A     &  $\Q$                 &  $-1/36$       \\
                     &   $(3,12)$          &     C     &      A     &  $\Q(\sqrt{3})$       &  $-23/1152$    \\
                     &   $(4,4)$           &     NC    &      A     &  $\Q$                 &  $-5/144$      \\
                     &   $(6,6)$           &     NC    &      A     &  $\Q$                 &  $-5/144$      \\
\hline
  $G_{29}$           &     $3$             &     NC    &     NA     &  $\Q(\sqrt{3})$       & $-323/12960$   \\
                     &     $4$             &     NC    &      A     &  $\Q$                 & $-13/160$      \\
\hline
  $G_{30}$           &     $3$             &     C     &      A     &  $\Q(\sqrt{5})$       & $-52/2025$     \\
                     &     $5$             &     C     &      A     &  $\Q(\sqrt{5})$       & $-41/1125$     \\
\hline
  $G_{31}$           &     $3$             &     NC    &      A     &  $\Q$                 & $-13/810$      \\
                     &     $5$             &     C     &      A     &  $\Q(\sqrt{5})$       & $-41/1125$     \\
\hline
\end{tabular}
\caption{Rough commensurability invariants and orbifold Euler
  characteristics, for CHL groups in
  $PU(3,1)$.} \label{fig:rough_invariants}
\end{table}

\end{document}